%% file: Intrinsic.tex
\documentclass[12pt,reqno]{amsart}
\usepackage{amssymb}
\usepackage{amsthm}
\usepackage{amsmath}
\usepackage{version}
\usepackage[margin=0.8in]{geometry}
\usepackage{graphicx}
\usepackage{enumerate}
\usepackage[usenames,dvipsnames] {color}

\def\m{\mathcal{M}}

\def\d{\mathbb{D}}
\def\c{\mathbb{C}}
\def\t{\mathbb{T}}
\def\r{\mathbb{R}}

\def\royal{\mathcal{R}}

\def\calr{\mathcal{R}}

\def\car#1{|#1|_{\rm car}}
\def\kob#1{|#1|_{\rm kob}}

\def\be{\begin{equation}}
\def\ee{\end{equation}}

\def\s0{s_0}
\def\p0{p_0}

\newcommand\nd{nondegenerate }

\DeclareMathOperator{\Car}{{\mathrm Car}}
\DeclareMathOperator{\Kob}{{\mathrm Kob}}
\DeclareMathOperator{\orb}{{\rm Orb}}

\numberwithin{equation}{section}
 \newtheorem{theorem}[equation]{Theorem}
 \newtheorem{corollary}[equation]{Corollary}

\newtheorem{lem}[equation]{Lemma}
\newtheorem{prop}[equation]{Proposition}
\newtheorem{thm}[equation]{Theorem}
\newtheorem{defin}[equation]{Definition}

\newtheorem{remark}[equation]{Remark}

\newtheorem{fact*}[equation]{Fact}

\DeclareMathOperator\aut{Aut}

\DeclareMathOperator{\cptwo}{{\mathrm CP}^1}
\DeclareMathOperator{\im}{\rm Im}

\newcommand\half{{\tfrac 12}}

\newcommand\id{{\mathrm{id}}}
\newcommand\idd{\mathrm{id}_\mathbb{D}}

\newcommand{\D}{\mathbb{D}}
\newcommand{\C}{\mathbb{C}}

\newcommand\fa{\mbox{ for all }}
\newcommand\fs{\mbox{ for some }}
\renewcommand\iff{\quad\Leftrightarrow\quad}

\def\set#1#2{\{ #1 \, : \, #2\}}
\def\id#1{{\rm id}_{#1}}

\newcommand{\inv}{^{-1}}

\newcommand{\ph}{\varphi}
\renewcommand\phi{\varphi}
\newcommand\al{\alpha}
\newcommand\ga{\gamma}
\newcommand\Ga{\Gamma}
\newcommand\de{\delta}

\newcommand\la{\lambda}

\newcommand\ups{\upsilon}
\newcommand\si{\sigma}
\newcommand\beq{\begin{equation}}
\newcommand\ds{\displaystyle}
\newcommand\eeq{\end{equation}}
\newcommand\df{\stackrel{\rm def}{=}}

\newcommand\red{\color{red}}

\newcommand\bbm{\begin{bmatrix}}
\newcommand\ebm{\end{bmatrix}}
\newcommand\bpm{\begin{pmatrix}}
\newcommand\epm{\end{pmatrix}}
\numberwithin{equation}{section}

\renewcommand\implies{\Rightarrow}
\makeindex

\begin{document}
\title[Intrinsic Directions in the Symmetrized Bidisc]{Intrinsic Directions, Orthogonality and Distinguished Geodesics in the Symmetrized Bidisc}
\author{Jim Agler}
\address{Department of Mathematics, University of California at San Diego, CA \textup{92103}, USA}
\email{jagler@ucsd.edu}
\thanks{Partially supported by National Science Foundation Grants
DMS 1361720 and 1665260, a Newcastle URC Visiting Professorship, the Engineering and Physical Sciences Research Council grant EP/N03242X/1 and London Mathematical Society Grant 41829} 

\author{Zinaida Lykova}
\address{School of Mathematics, Statistics and Physics, Newcastle University, Newcastle upon Tyne
 NE\textup{1} \textup{7}RU, U.K.}
\email{Zinaida.Lykova@ncl.ac.uk}

\author{N. J.  Young}
\address{School of Mathematics, Statistics and Physics, Newcastle University, Newcastle upon Tyne NE1 7RU, U.K.
{\em and} School of Mathematics, Leeds University,  Leeds LS2 9JT, U.K.}
\email{Nicholas.Young@ncl.ac.uk}
\date{6th December 2020}


\begin{abstract} 
The symmetrized bidisc 
\[
G \stackrel{\rm{def}}{=}\{(z+w,zw):|z|<1,\ |w|<1\},
\]
under the Carath\'eodory metric, is a complex Finsler space of cohomogeneity $1$  
in which the geodesics, both real and complex, enjoy a rich geometry.  As a Finsler manifold, $G$ does not admit a natural notion of angle, but we nevertheless show that there {\em is} a notion of orthogonality.  The complex tangent bundle $TG$ splits naturally into the direct sum of two line bundles, which we call the {\em sharp} and {\em flat} bundles, and which are geometrically defined and therefore covariant under automorphisms of $G$.   Through every point of $G$ there is a unique complex geodesic of $G$ in the flat direction, having the form
\[
F^\beta  \stackrel{\rm{def}}{=}\{(\beta+\bar\beta z,z)\ : z\in\D\}
\]
for some $\beta \in\D$, and called a {\em flat geodesic}.
We say that a complex geodesic \emph{$D$ is orthogonal} to a flat geodesic $F$ if $D$ meets $F$ at a point $\lambda$ and the complex tangent space $T_\lambda D$ at $\lambda$ is in the sharp direction at $\lambda$. We prove that a geodesic $D$ has the closest point property with respect to a flat geodesic $F$ if and only if $D$ is orthogonal to $F$ in the above sense.
Moreover, $G$ is foliated by the geodesics in $G$ that are orthogonal to a fixed flat geodesic $F$.

\end{abstract} 

\maketitle
\tableofcontents
\input intro

\input Car-Kob

\input geodesics

\input purebalgeos

\input orthogonality

\input distinguished

\input closestpoint

\input Intrinsic.bbl
\end{document}

%% file: intro.tex


\section{Introduction}\label{intro}

This paper concerns {\em geodesics} in the symmetrized bidisc 
$G$, the domain in $\C^2$ defined in the abstract. There are two different notions of geodesic, to wit, {\em complex} and {\em real}.
\index{geodesic!complex}
\index{geodesic!real}
A {\em complex geodesic} in $G$  is defined to be the range of a holomorphic map $h:\D\to G$ which is isometric with respect to the Poincar\'e metric on the unit disc $\D =\set{z \in \c}{|z|<1}$ and the hyperbolic metric on $G$.  It is a fact that such an $h$ is an isometry if and only if  $h$ has a holomorphic left inverse (see, for example, \cite{jp04}).  {\em Real geodesics} are paths which locally minimize lengths, in the usual sense of real differential geometry or metric geometry (see \cite[Definition 2.5.27]{BBI2001}).  When we speak of a `geodesic' in $G$, without qualification, we shall mean a {\em complex} geodesic.

We study the geometry of $G$ as a hyperbolic complex space in the sense of Kobayashi \cite{kob98}.  It is well known that on $G$ the {\em Carath\'eodory pseudometric} $\car{\delta}$ and {\em Kobayashi pseudometric} $\kob{\delta}$, where $\delta$ is a datum{\footnote{That is, either a pair of points in $G$ or an element of the complex tangent bundle of $G$} in $G$, are in fact metrics, and moreover they coincide  \cite{ay2004,jp04}. Thus there is a natural metric $d_G$   on $G$, given, for a datum $\delta= (\lambda, \mu)$ in $G$, by
\[
d_G(\delta) = \kob{\delta}=\car{\delta}.
\]
We call $d_G$  the {\em hyperbolic metric} on $G$. 
\index{hyperbolic metric on $G$}

A fundamental fact of the complex geometry of $G$ is that, for every point $\lambda\in G$ and every one-dimensional subspace  $\m \subseteq \c^2$, there exists a {\em unique} complex geodesic $D$ in $G$
such that $\lambda \in D$ and $T_\lambda D=\m$,
where $T_\lambda D$ denotes the complex tangent space to $D$ at the point $\lambda$.
Thus there is a natural one-to-one correspondence between the geodesics in $G$ passing through $\lambda$ and the points of  $\cptwo$, the {\em complex projective line.} We regard $\cptwo$ as the set of one-dimensional subspaces of $\c^2$ and refer to its elements as {\em directions}. 
\index{$\cptwo$}
\index{direction}
We shall exploit this correspondence to relate the local and global hyperbolic geometries of $G$.

As there are qualitative differences among the complex geodesics of $G$, so are there differences among the directions at a point. 
To illustrate this anisotropy we first note that 
 there exists a unique complex geodesic $\calr$ in $G$  that is invariant under all automorphisms of $G$.
Referred to as the \emph{royal} variety, $\royal$ is defined to be the set
\index{royal!variety}  \index{$\royal$} \index{$\calr$}
\[
 \royal = \{(2z, z^2), z\in \d \} = \{(s,p)\in G: s^2=4p\}.
\]
Moreover,
for every $\la\in\calr$, there exists a unique geodesic $F_\la$ in $G$  having a nontrivial stabilizer in $\aut G$, the group of biholomorphic self-maps of $G$, and such that
\[
F_\la\cap\calr = \{\la\}.
\]
Thus the set $\{F_\la: \la \in \calr \}$ is a hyperbolically identifiable class of geodesics. We call them \emph{flat} geodesics. They have the form
\be
F^\beta \df \{(\beta+\bar\beta z,z): z\in\d\}
\ee
for some $\beta \in\d$.  
 Furthermore through each point $\lambda \in G$ there passes a unique flat geodesic.
Consequently, at each point $\lambda \in G$, there exists a unique direction
  $\flat_\lambda \in \cptwo$ with the property that
\[
T_\lambda F = \flat_\lambda,
\]
where $F$ is the flat geodesic passing through $\lambda$. We refer to $\flat_\lambda$ as the \emph{flat direction at $\lambda$}.
\index{flat direction}
\index{$\flat_\lambda$}
 The flat direction in $G$ is covariant under automorphisms of $G$,  
\index{covariant}
in the sense that for every $\lambda \in G$ and $\gamma \in \aut G$,
\[
  \flat_{\gamma(\lambda)}=\gamma'(\lambda)\flat_\lambda.
\]

Another covariant direction is the {\em sharp} direction \cite[Subsection 1.4]{aly2017}. 
\index{sharp!direction}
For a point $\lambda \in G$, let $\orb(\lambda)$ denote the orbit of $\lambda$ under the action of $\aut G$.
In \cite[Theorem 1.6]{aly2017} the authors showed that if $\lambda \in \royal$, then $ \orb(\lambda) = \royal$,
and if  $\lambda \in G \setminus\royal$, then
$\orb(\lambda)$  is a smooth properly embedded 3-dimensional real manifold in $G$. Noting that, in both cases, $T_\la\orb(\lambda)$  contains a unique one-dimensional complex subspace $\m$, we may define an intrinsic direction $\sharp_\lambda$, the \emph{sharp direction at $\lambda$}, by $\sharp_\lambda \df \m$.
\index{$\sharp_\lambda$}
Observe that, for $\lambda \in \royal$, since $\orb(\la)=\calr$,
\be\label{shroyal}
\sharp_\lambda = T_\la\orb(\la) =T_{\lambda}\royal.
\ee
In this paper, besides the sharp and flat directions we shall encounter other special directions at a point. We now describe the results of the paper.

If $F$ is a flat geodesic in $G$ and $D$ is a general geodesic in $G$, we say that \emph{$D$ is orthogonal to $F$}
\index{orthogonal}
 if $D$ meets $F$ at a point $\lambda$ and $D$ points in the sharp direction at $\lambda$, that is, $T_\lambda D =\sharp_\lambda$. The motivation for this terminology is  the following result (Corollary \ref{cpp.cor.20} below). 
\begin{thm}\label{int.thm.10}
 If $F$ is a flat geodesic in $G$ and $D$ is an arbitrary geodesic in $G$, then $D$ is orthogonal to $F$ if and only if $D$ meets $F$ at a point $\lambda_0$ and, for some point $\mu \in D\setminus \{\lambda_0\} $ (equivalently, for every point $\mu \in D\setminus \{\lambda_0\}$), 
\[
d_G(\lambda_0,\mu)=\inf_{\lambda\in F} d_G(\lambda,\mu).
\]
\end{thm}

The royal geodesic meets any flat geodesic in a unique point. Therefore, by the above definition and equation \eqref{shroyal}, $\royal$ is orthogonal to every flat geodesic. If $D\not= \royal$, then for $D$ to be orthogonal to a flat geodesic, it must be the case that $D$ is of a very special type which we now describe (see Corollary \ref{sharp.cor.20} and \cite[Lemma 9.8]{aly2016}). We say that a geodesic $D$ is \emph{purely balanced}
\index{geodesic!purely balanced}
 if $D^- \cap \royal^-$ consists of exactly two points $\xi_1$ and $\xi_2$ both of which lie in 
$\partial \royal= \{(2z,z^2): z \in \t \}$\footnote{The closure and boundary symbols relate to $D$ and $\calr$ as subsets of $\c^2$.}.
Here  $\t$ denote the unit circle $\set{z \in \c}{|z|=1}$. The points $\xi_1$ and $\xi_2$ are referred to as the \emph{royal points of $D$}. 
\index{royal!points}
In Section \ref{purebalgeos} of this paper we  describe the purely balanced geodesics using properties of hyperbolic automorphisms of the disc.\footnote{$m\in \aut \d$ is said to be \emph{hyperbolic} if $m$ has two distinct fixed points on $\t$.} 
\index{hyperbolic automorphism}
This description is used to show that the purely balanced directions at a point in $G \setminus \royal$ form a simple smooth curve in $\cptwo$ connecting two `exceptional' directions and containing the sharp direction as its midpoint (cf. Corollary \ref{hyp.cor.10}).

We further exploit the description of the purely balanced geodesics to prove the following result in Section \ref{orthogonality} of the paper (Theorem \ref{sharp.thm.10}).

\begin{thm}\label{int.thm.20}
 If $F$ is a flat geodesic in $G$, then $G$ is foliated by the geodesics in $G$ that are orthogonal to $F$.
\end{thm}

As a corollary of Theorems \ref{int.thm.10} and \ref{int.thm.20} we obtain in  Corollary \ref{cpp.cor.10} the following result.
\begin{corollary}\label{int.cor.10}
If $F$ is a flat geodesic in $G$ and $\mu\in G$, then there exists a unique point $\lambda_0 \in F$ such that
\[
d_G(\lambda_0,\mu)=\inf_{\lambda\in F} d_G(\lambda,\mu).
\]
\end{corollary}
 \cite[Lemma 9.8]{aly2016}  contains the statement that the shortest distance from $\mu$ to $F$ is always attained; we now see that it is attained at a unique point of $F$.

In Section \ref{distinguished} we turn the tables -- we fix a  purely balanced geodesic $D$ and study the flat geodesics $F$ such that $D$ is orthogonal to $F$.
 If $D$ is a geodesic, we let $\sharp(D)$ denote the set of $\lambda$ in $D$ such that $D$ is orthogonal to the flat geodesic passing through $\lambda$ (that is, $T_\lambda D = \sharp_\lambda$.) We say that a curve $\Xi$ is a \emph{simple real geodesic in $G$} if $\Xi = h(C)$, where $h:\d \to G$ is an isometry when $\d$ and $G$ are equipped with their respective hyperbolic metrics, and $C$ is a real hyperbolic 
geodesic\footnote{that is, $C$ is either a line segment or a circular arc in $\d$ that intersects $\t$ orthogonally.} in $\d$.
Each purely balanced geodesic $D$ contains a unique real geodesic $\Xi_D$ such that the endpoints of $\Xi_D$ are the royal points of $D$ (cf. Proposition \ref{dis.prop.40}).
\begin{thm}\label{int.thm.30}
 If $D$ is a purely balanced geodesic in $G$, then $\sharp(D) = \Xi_D$, that is, the set of $\lambda \in D$ such that $T_{\lambda}D = \sharp_{\lambda}$ is the simple real geodesic in $G$ whose endpoints are the two royal points of $D$.
\end{thm}

As a corollary of Theorem \ref{int.thm.30} we obtain in Section \ref{distinguished} the following result.
\begin{thm}\label{int.thm.40}
{\em (cf. Theorem \ref{dis.thm.20})} $G \setminus \royal$ is foliated by the simple real geodesics in $G$ whose endpoints lie in $\partial \royal$.
\end{thm}

Another special direction at a point $P$ of $G$, this time in the {\em real} tangent bundle of $G$, we call the {\em distinguished} direction.  It 
 is tangent to a special real geodesic $\gamma$ through $P$ with the property that, for every pair of points on $\gamma$, there are at least two inequivalent solutions of the corresponding Carath\'eodory extremal problem.
These distinguished geodesics foliate $G\setminus \mathcal{R}$, where $\mathcal{R}$ is the `royal variety' $\{(2z,z^2): |z|<1\}$.

The theory of  the symmetrized bidisc, and cognate domains like the tetrablock, has been extensively developed over the last 20 years by numerous authors. We shall require some results from this theory,  
many of which can be found in  \cite{jp} and  \cite[Appendix A]{aly2016}.
These domains have a rich complex geometry and function theory, as well as applications to operator theory: see, besides many other papers, \cite{KZ1,KZ2,jp,ALY18,KMcC,ballsau,bhatta,bhatta17}.

%% file: Car-Kob.tex
\section{The hyperbolic metric on a Lempert domain}\label{CandK}

In this section we describe our terminology for the Carath\'eodory and Kobayashi extremal problems on  an open set $U$ in $\c^n$ and introduce the hyperbolic metric on a Lempert domain.

For open sets $U \subseteq \c^{n_1}$ and $V \subseteq \c^{n_2}$ we denote by $V(U)$ the set of holomorphic mappings from $U$ into $V$.

If $U$ is an open set in $\c^n$, then by a \emph{datum in} $U$
 we mean an ordered  pair $\delta$ where either $\delta$ is \emph{discrete}, that is, has the form
\[
\delta =(s_1,s_2)
\]
where $s_1,s_2 \in U$, or $\delta$ is \emph{infinitesimal}, that is, has the form
\[
\delta = (s,v)
\]
where $s \in U$ and $v\in\c^n$. 

If $\delta$ is a datum, we say that $\delta$ is \emph{degenerate} if either $\delta$ is discrete and $s_1=s_2$ or $\delta$ is infinitesimal and $v=0$. Otherwise, we say that $\delta$ is \emph{nondegenerate}.

An infinitesimal datum in $U$ is the same thing as a point of the complex tangent bundle $TU$ of $U$.

For $F\in \Omega(U)$, $s\in U$, and $v \in \c^{n_1}$, the directional derivative $D_v F(s) \in \c^{n_2}$ is defined by
\index{$D_vF$}
\[
D_v F(s) = \lim_{z \to 0} \frac{F(s+zv) - F(s)}{z}.
\]

 If $U$ and $\Omega$ are domains, $F\in \Omega(U)$, and $\delta$ is a datum in $U$, we define a datum $F(\delta)$ in $\Omega$ by
 \[
 F(\delta) = (F(s_1),F(s_2))
 \]
 when $\delta=(s_1,s_2)$ is discrete and by
 \[
 F(\delta) =(F(s),D_v F (s))
 \]
 when $\de=(s,v)$ is infinitesimal. 

 For any datum $\delta$ in $\d$, we define $|\delta|$ to be the Poincar\'e distance or metric at $\de$ in the discrete or infinitesimal case respectably, that is
\[
|\delta| = \tanh\inv \left|\frac{z_1 -z_2}{1-\bar{z_2}z_1}\right|
\]
when $\delta=(z_1,z_2)$ is discrete\footnote{In \cite[Chapter 3]{aly2016} we used a different notation in that we omitted $\tanh\inv$; this makes no essential difference in the present context, but ensures that the Carath\'eodory pseudodistance is the inner pseudodistance determined by the Carath\'eodory pseudometric \cite{jp}, and similarly for the Kobayashi pseudodistance.},
and by
\[
|\delta|=\frac{|v|}{1-|z|^2}
\]
when $\delta = (z,v)$ is infinitesimal.
See \cite{kr} for terminology and theory in several complex variables.
\\ 

{\bf  The Carath\'eodory extremal problem.}
\index{Carath\'eodory!extremal!problem}
{\em For a domain  $U$ in $\c^n$ and a nondegenerate datum $\delta$  in $U$, compute the quantity $\car{\delta}$
defined by}
\beq\label{defCmetric}
\car{\delta}=\sup_{F\in \d(U)} |F(\delta)|.
\eeq
We shall refer to this problem as $\Car\delta$ 
and will say that \emph{$C$ solves} $\Car \delta$ if $C\in \d(U)$ and
\[
\car{\delta} = |C(\delta)|.
\]

$\car{\cdot}$ is a pseudometric on $U$, called the {\em Carath\'eodory pseudometric} in the case of infinitesimal $\de$ and the  {\em Carath\'eodory pseudodistance } in the case of discrete datums $\de$.
\index{ Carath\'eodory!pseudometric}
\index{ Carath\'eodory!pseudodistance}

It is easy to see with the aid of Montel's theorem that, for every \nd datum $\de$ in $U$, there does exist $C\in\d(U)$ which solves $\Car\de$.  Such a $C$ is called a {\em Carath\'eodory extremal function } for $\de$.
\index{Carath\'eodory!extremal!function}

\begin{defin}\label{predef20}
We say that a domain $U$ in $\c^n$ is \emph{weakly hyperbolic}
\index{weakly hyperbolic}
 if $\car{\delta} >0$ for every nondegenerate datum $\delta$ in $U$. Equivalently, for every nondegenerate datum $\delta$ in $U$, there exists a bounded holomorphic function $F$ on $U$ such that $F(\delta)$ is a nondegenerate datum in $\c$.
\end{defin}

\noindent{\bf  The Kobayashi extremal problem.} 
\index{Kobayashi!extremal!problem}
{\em For  a domain $U$ in $\c^n$ and a nondegenerate datum $\delta$ in $U$, compute the quantity $\kob{\delta}$
 defined by}
\beq\label{defkob2}
\kob{\delta}=\inf_{\substack{f\in U(\d)\\ f(\zeta) =\delta}} |\zeta|.
\eeq
We shall refer to this problem as $\Kob \delta$ 
and will say that \emph{$k$ solves} $\Kob\delta$
 if $k\in U(\d)$ and there exists a datum $\zeta$ in $\d$ such that $k(\zeta) = \delta$ and
\[
\kob{\delta} = |\zeta|.
\]

On infinitesimal datums $\kob{\cdot}$ is a pseudometric, called the {\em Kobayashi pseudometric} or the {\em Kobayashi-Royden pseudometric} \cite[Chapter 3]{jp}.
\index{Kobayashi!pseudometric}
The quantity $\kob{\cdot}$ is not necessarily a pseudodistance on discrete datums (it can fail to satisfy the triangle inequality).  The {\em Kobayashi pseudodistance} on $U$
\index{Kobayashi!pseudodistance}
is defined to be the largest pseudodistance on $U$ majorized by  $\kob{\cdot}$.

Note that the infimum in the definition \eqref{defkob2} of $\kob{\de}$ is attained if $U$ is a taut domain, where $U$ is said to be {\em taut} if $U(\d)$ is a normal family.   In particular $\kob{\de}$ is attained when $U=G$  \cite[Section 3.2]{jp}.  Any function which solves $\Kob\de$ is called a {\em Kobayashi extremal function for $\de$}.

The Kobayashi and  Carath\'eodory pseudometrics are {\em invariant}, that is, they are invariant under automorphisms of $G$; see \cite{jp} for an up-to-date account of such pseudometrics.

Let $U$ be a domain in $\c^n$ and $\delta$ a nondegenerate datum in $U$. The solutions to $\Car \delta$ and $\Kob\delta$ are never unique, for if $m$ is a M\"obius transformation of $\d$, then $m\circ C$ solves $\Car\delta$ whenever $C$ solves $\Car\delta$ and $f \circ m$ solves $\Kob\delta$ whenever $f $ solves $\Kob\delta$. This suggests the following definition. 
 
\begin{defin}\label{predef30}
Let $U$ be a domain in $\c^n$ and let $\delta$ be a nondegenerate datum in $U$. We say that \emph{the solution to $\Car\delta$ is essentially unique}, if whenever $F_1$ and $F_2$ solve $\Car\delta$ there exists a M\"obius transformation $m$ of $\d$ such that $F_2=m \circ F_1$.
We say that \emph{the solution to $\Kob\delta$ is essentially unique} if the infimum in equation \eqref{defkob2} is attained and, 
 whenever $f_1$ and $f_2$ solve $\Car\delta$ there exists a M\"obius transformation $m$ of $\d$ such that $f_2=f_1 \circ m$.
\end{defin}

In honor of Lempert's seminal theorem \cite{lem81} we adopt the following definition.
\begin{defin}\label{predef50} 
A domain $U$ in $\c^n$ is a \emph{Lempert domain}
\index{Lempert!domain}
 if 
\begin{enumerate}[\rm (1)]
\item $U$ is weakly hyperbolic,
\item  $U$ is taut and 
\item $\car{\delta}=\kob{\delta}$ for every nondegenerate datum $\delta$ in $U$.
\end{enumerate}
\end{defin}
Thus,  for  a Lempert domain $U$ in $\c^n$, there is a natural metric $d_U$   on $U$, given, for a nondegenerate infinitesimal datum $\delta$ in $U$,
by
\[
d_U(\delta) = \kob{\delta}=\car{\delta}.
\]
The metric $d_U$ is called the {\em hyperbolic metric} on $U$.
\index{hyperbolic metric on a Lempert domain}\\

\noindent{\bf The hyperbolic metric on the symmetrized bidisc.}
The {\em symmetrization map} $$\pi:\c^2\to \c^2$$  is defined by
\[
\pi\big((z_1,z_2)\big) = (z_1+z_2,z_1z_2) \qquad \fa z_1,z_2\in\c.
\]
Thus  the symmetrized bidisc  $G$ is $ \pi(\d^2)$.

We adopt the co-ordinates $s=z_1 + z_2$ and $p=z_1z_2$. For $s,p \in \c$, 
\be\label{schurcrit}
(s,p) \in G \mbox{ if and only if } |s-\bar s p| < 1-|p|^2
\ee
(for example \cite[Lemma 7.1.3]{jp}).

The symmetrized bidisc $G$ is a Lempert domain in $\c^2$ with the hyperbolic metric $d_G$  \cite{ay2004}.

%% file: geodesics.tex

\section{Automorphisms, complex geodesics and directions in $G$}\label{geodesics}

\subsection{The automorphism group of $G$}

For a domain  $\Omega$ in $\C^d$, an {\em automorphism} of $\Omega$ is defined to be a 
 biholomorphic self-map of $\Omega$. The group of all 
 automorphisms  of $\Omega$ under composition will be denoted by  $\aut \Omega$.
Here is a description of $\aut G$ in terms of $\aut \d$.

\begin{prop}\label{pre.prop.5}
For every $b\in \aut \d$, there is a unique automorphism $\gamma_b$ of $G$ satisfying
\be\label{pre.5}
\gamma_b (\pi (z_1,z_2)) = \pi(b(z_1),b(z_2)) \qquad  \fa z_1,z_2 \in \d.
\ee
Furthermore, the map $f: \aut \d \to \aut G$ given by 
\[
f( b ) = \gamma_b \; \qquad  \fa b \in \aut \d
\]
is a continuous isomorphism with respect to the compact-open topologies on $\aut \d$ and $\aut G$. 
\end{prop} 
A proof can be found in  \cite[Section 7.1]{jp} or  \cite[Theorem 4.1]{AY08}.

The following corollary of Proposition \ref{pre.prop.5} is an immediate consequence of the fact that every automorphism of $\d$ extends to be analytic in a neighborhood of $\d^-$.
\begin{corollary}\label{extends}

\begin{enumerate}[\rm (1)]
\item Every automorphism $\ga$ of $G$ extends in a unique way to a self-homeomorphism $\tilde\ga$ of the closure of $G$ in $\c^2$; 
\item $\tilde\ga(\calr^-)=\calr^-$;
\item the classes of flat and purely balanced geodesics in $G$ are invariant under the automorphisms of $G$.
\end{enumerate}
\end{corollary}

\subsection{The complex geodesics in $G$}
By a {\em complex geodesic} in $G$ we mean a subset $D \subset G$ such that there exists a hyperbolic isometry $h:\d \to G$ such that  $D=h(\d)$. It is known that through any two points of $G$ there passes a unique geodesic (see \cite[Theorem 0.3]{AY06}). We note that $h:\d \to G$ is a hyperbolic isometry if and only if $h$ is holomorphic and there exists a holomorphic map 
$\Phi: G \to \d $ such that\footnote{Here $\id{\d}$ denotes the automorphism defined by $\id{\d}(z) = z,\ \ z \in \d$.}  $\Phi \circ h=\id{\d}$.  
Moreover, when $D=h(\d)$ is a geodesic, $h$ is rational and extends to be holomorphic on a neighborhood of $\d^-$, so that, in particular, $D^- =h(\d^-)$.

There are qualitative differences between the geodesics of $G$. Indeed, in \cite[Chapter 7]{aly2016} five distinct types of geodesic are identified, namely \emph{royal, flat, purely unbalanced, purely balanced,} and \emph{exceptional}. Each of the types can be characterized in terms of qualitative properties of the Carath\'eodory extremal problem associated with the geodesic, and in a variety of purely geometric ways as well.

Here, we adopt a description of the five types of geodesics in terms of the geometry of the distinguished boundary\footnote{That is, the smallest closed  subset ${\mathrm M}$ of $G^-$ with the property that every holomorphic function defined on a neighborhood of $G^-$ attains its maximum over $G^-$ at a point in ${\mathrm M}$.} of $G$,
\index{distinguished boundary}
\index{${\mathrm M}$}
which can be shown to be the set  in $\c^2$ defined by
\[
{\mathrm M} = \pi (\t^2).
\]
Topologically, ${\mathrm M}$ is a M\"obius band. The {\em edge} ${\mathrm E}$ of the M\"obius band is the set in $\c^2$ 
\[
{\mathrm E} = \pi(\set{(z,z)}{z \in \t}).
\]
\index{$\mathrm E$}
\index{edge}
Note that ${\mathrm E} =\partial \royal$.

\subsection{The five types of geodesic}\label{5types}
We summarize the definitions of the types of geodesic introduced in \cite{aly2016} and state a few of their properties.  In that paper \cite{aly2016} the types were defined in Subsection 3.3, but here we shall rather use a geometric characterization of them given in \cite[Chapter 7]{aly2016}.
 \begin{enumerate}
 \item The \emph{royal} geodesic is the set
\index{geodesic!royal}
\index{royal!geodesic}
 \[
 \royal = \pi (\set{(z,z)}{z\in \d}).
 \]
 It is the only geodesic left invariant by every automorphism in $\aut G$. Furthermore\footnote{Indeed, $\royal^-$ is the polynomially convex hull of ${\mathrm E}$.},
 \[
 \royal^- \cap {\mathrm M} = {\mathrm E}.
 \]
 Note also that $ \royal^-$ is the disjoint union of $\royal$ and ${\mathrm E}$.

If $D$ is a geodesic we define the {\em royal points} of $D$ to be the elements of the set $D^- \cap \royal^-$.
\index{royal!points}
Since $\royal^- = \royal \cup {\mathrm E}$, each royal point must either be in $\royal$ or in $ {\mathrm E}$.

 \item  A geodesic $D$ is said to be \emph{flat}
\index{geodesic!flat}
 if  $D^- \cap \royal^-$ consists of a single point, and that point lies\footnote{Equivalently, a geodesic $D$ is flat if and only if $D^- \cap {\mathrm E} = \emptyset$.} in $\royal$.
Thus, a geodesic $D$ is flat if it has a unique royal point and that point lies in $\royal$ (rather than $\mathrm E$).
The nomenclature  ``flat geodesic" reflects the fact that a
 geodesic  $D$ is flat if and only if $D$ lies in a set of the form $\lambda+ \m$ where $\lambda \in G$ and $\m$ is a one-dimensional complex subspace of $\c^2$. There exists a unique flat geodesic passing through each point of $G$, or in other words, the flat geodesics foliate $G$.
     The royal geodesic is the only geodesic in $G$ that meets every flat geodesic.
          \item A geodesic $D$ is said to be \emph{purely balanced} if $D^- \cap \royal^-$ consists of exactly two points both of which lie in ${\mathrm E}$,
that is, $D$ has exactly two royal points, and they lie in ${\mathrm E}$. We shall give a concrete formula for the general purely balanced geodesic using hyperbolic automorphisms of $\d$ in Section \ref{purebalgeos}.
    \item A geodesic $D$ is said to be \emph{purely unbalanced} 
\index{geodesic!purely unbalanced}
if $D^- \cap \royal^-$ consists of exactly two points, one of which lies in $\royal$ and one of which lies in ${\mathrm E}$, that is, $D$ has exactly two royal points, one of which lies in $\royal$ and one of which lies in ${\mathrm E}$. Generically, the geodesic that passes through two distinct points in $G$ is purely unbalanced.
 
          \item A geodesic $D$ is said to be \emph{exceptional}
\index{geodesic!exceptional}
 if $D^- \cap \royal^-$ consists of exactly one point and that point lies in ${\mathrm E}$, that is, $D$ has exactly one royal point, and that point lies in ${\mathrm E}$.
 \end{enumerate}

\subsection{Directions}
We shall denote a direction at a point in $G$ (that is, an element of $\cptwo$) by 
\[
v\c =\set{zv}{z \in \c}\qquad \mbox{ where } v\in \c^2 \setminus\{0\}.
\]
 The following pleasing state of affairs is central to this paper  (see \cite[Theorem 4.6]{aly2016}).
\begin{thm}\label{pre.thm.10}
 For every  point $\lambda$ in $G$ and every nonzero vector $v$ in $\c^2$, there exists a unique geodesic $D^{v}$ in G such that
\[
\lambda \in D^v\ \text{ and }\ v \in T_\la D^{v}.
\]
Furthermore, if $v$ and $w$ are any two nonzero vectors in $\c^2$,
$D^{v}=D^{w}$ if and only if $v\c =w\c$. In particular,
\be\label{pre.10}
 v\c \mapsto D^{v}
\ee
is a well-defined bijection between directions in $\cptwo$ and geodesics in $G$ that pass through $\lambda$.
\end{thm}
Evidently, Theorem \ref{pre.thm.10} implies that just as there are five qualitatively distinct types of geodesic, so are there five qualitatively distinct types of direction at a point. Let us say that a direction $v\c$ is of a particular type if the geodesic $D^v$ that corresponds to $v\c$ via the map \eqref{pre.10} is of that same type.
\subsubsection{The flat direction}
As noted above, through each point of $G$ there passes a unique flat geodesic in $G$. Hence, for each $\lambda\in G$ there is a unique direction $\flat_\lambda \in \cptwo$  defined by the following procedure.
\begin{enumerate}
\item Let $D$ be the unique flat geodesic in $G$ passing through $\lambda$.
\item Choose any nonzero vector $v$ in $T_\lambda D$.
\item Let $\flat_\lambda = v\c$.
\end{enumerate}
We refer to $\flat_\lambda$ as the \emph{flat direction at $\lambda$} and say that a nonzero vector $v\in \c^2$ \emph{points in the flat direction at $\lambda$} if $v \in \flat_\lambda$.
\subsubsection{The sharp direction}
\index{sharp!direction}
In \cite[Theorem 1.6]{aly2017}, the authors showed that if $\lambda \in G$ and $\orb(\lambda)$ 
\index{$\orb(\la)$}
denotes the orbit of $\lambda$ under $\aut G$, then
\begin{enumerate}
\item $\orb(\lambda) = \royal$ if $\lambda \in \royal$, and
\item $\orb(\lambda)$ is a smooth properly embedded 3-dimensional real manifold in $G$ if $\lambda \not\in \royal$.
\end{enumerate}

As a consequence of this result, for each $\lambda\in G$, we may define $\sharp_\lambda \in \mathrm{CP^1},$ by the following procedure.
\begin{enumerate}
\item If $\lambda \in \royal$, then $\sharp_\lambda = v\c$, where $v$ is any nonzero vector in $T_\la\royal$.
    \item Otherwise, if $\lambda \not\in \royal$, let $\sharp_\lambda = v\c$, where $v$ is any vector such that both $v$ and $iv$ are in the real tangent space to $\orb(\lambda)$ at $\lambda$.
\end{enumerate}
We refer to $\sharp_\lambda$ as the \emph{sharp direction at $\lambda$} and we say that a nonzero vector $v\in \c^2$ \emph{points in the sharp direction at $\lambda$} if $v \in \sharp_\lambda$. Evidently, with this language,  if $\lambda\in G$ and $v$ is a nonzero vector in $\c^2$, then $v$ points in the sharp direction at $\lambda$ if and only if $v$ is in the unique nonzero complex subspace of the real tangent space to $\orb(\lambda)$ at $\lambda$.

\subsubsection{Other directions}
For each point $\lambda$ of $G$ the previous two subsections constructed  a unique pair of directions $\flat_\lambda$ and $\sharp_\lambda$. These directions, being defined geometrically,  are covariant in the sense that if $\gamma$ is an automorphism of $G$, then
\[
\flat_{\gamma(\lambda)} = \gamma'(\lambda)\ \flat_\lambda \ \ \text{ and }\ \ \sharp_{\gamma(\lambda)} = \gamma'(\lambda)\ \sharp_\lambda \quad \fa \lambda\in G.
\]
The consideration of covariant directions other than the flat and sharp directions is complicated by issues of both existence and uniqueness. At a point $\lambda \in G$ there may not be any directions of a given type. Alternatively, there may be multiple directions of a given type. We briefly savour the low-hanging fruit regarding these issues.
\begin{enumerate}
\item As the royal geodesic $\royal$ is unique, if $\lambda \in G$, then there exists a royal geodesic passing through $\lambda$ if and only if $\lambda\in \royal$. Therefore, there exists a `royal direction' at $\lambda$ only when $\lambda\in \royal$. However, in that case, $T\royal_\lambda= \sharp_\lambda$, that is, a `royal direction' exists only if it is the sharp direction. For this reason we shall not henceforth use the term royal direction.
    \item According to our definition, a direction $v\c$ is purely balanced at $\lambda\in G$ if the unique geodesic $D$ passing through $\lambda$ satisfying $T_\lambda D =v\c$ is purely balanced.  Notice from item (4) above in the discussion of the five types of geodesic that no purely balanced geodesic can meet $\royal$. Therefore, if $\lambda \in \royal$, then there are no purely balanced directions at $\lambda$.
        \item On the other hand, if $\lambda \in G\setminus \royal$, it turns out that  there is a one-parameter family of purely balanced directions at $\lambda$, see (7) below.
                \item A direction $v\c$ is \emph{exceptional at $\lambda\in G$} if the unique geodesic $D$ passing through $\lambda$ satisfying $T_\lambda D =v\c$ is exceptional.  Notice from item (5) above in the discussion of the five types of geodesics that no exceptional geodesic can meet $\royal$. Therefore, if $\lambda \in \royal$, then there are no exceptional directions at $\lambda$.
        \item On the other hand, if $\lambda \in G\setminus \royal$, it turns out that there are exactly two exceptional\footnote{The musical symbol $\natural$ is read `natural'.}
 directions $\natural^1_\lambda$ and $\natural^2_\lambda$ at $\lambda$. We let $\natural_\lambda$ denote the set of these two directions, that is, $\natural_\lambda=\{\natural^1_\lambda,\natural^2_\lambda\}$, see (7) below.
            \item Strictly speaking, the exceptional directions are not covariant. However, the {\em pair} of exceptional direction is covariant: if $\lambda \in G \setminus \royal$ and $\gamma \in \aut G$, then
                 \[
                 \natural_{\gamma(\lambda)} = \gamma'(\lambda)(\natural_\lambda).
                 \]
              \item If $\lambda \in G \setminus \royal$ then the purely balanced directions at $\lambda$ form a simple smooth curve connecting $\natural^1_\lambda$ and $\natural^2_\lambda$ in $\cptwo$ (cf. Corollary \ref{hyp.cor.10}).
                  \item If $\lambda\in G$, we say a direction $v\c$ is \emph{purely unbalanced} at $\lambda$ if the unique geodesic $D$ passing through $\lambda$ satisfying $T_\lambda D=v\c$ is purely unbalanced.
                      \item If $\lambda\in \royal$ and $v\c$ is a direction at $\lambda$, then either $v\c=\flat_\lambda$, $v\c=\sharp_\lambda$, or $v\c$ is purely unbalanced.
                      \item If $\lambda \in G \setminus \royal$ and $v\c$ is a direction at $\lambda$, then either $v\c=\flat_\lambda$, $v\c \in \natural_\lambda $, $v\c$ is purely balanced, or $v\c$ is purely unbalanced or  exceptional. The sharp direction is the midpoint of the purely balanced curve alluded to in item (7) above (cf. Corollary \ref{sharp.cor.40}).
\end{enumerate}

%% file: purebalgeos.tex

\section{Hyperbolic automorphisms in $\d$ and purely balanced geodesics in $G$}\label{purebalgeos}

 \subsection{Hyperbolic automorphisms of $\d$}
 If $\alpha\in \d$ we define $b_\alpha \in \aut \d$ by
\index{$b_\alpha$}
 \[
 b_\alpha(z) = \frac{z-\alpha}{1-\bar\alpha z},\qquad z \in \d,
 \]
 and when $\tau \in \t$, we define $r_\tau \in \aut \d$ by
\index{$r_\tau$}
 \[
 r_\tau (z)= \tau z, \qquad z \in \d.
 \]
Each $m\in \aut \d$ can be uniquely represented in the form
 \be\label{hyp.5}
 m = r_\tau \circ b_\alpha
 \ee
for some $\alpha \in \d$ and $\tau \in \t$.

 If $m \in \aut \d\setminus \{\id{\d}\}$
  then $m$ can be viewed as a continuous self-map of $\d^-$ and it is well known that exactly one of the following cases occurs.
 \begin{enumerate}[\rm (i)]
 \item $m$ is \emph{elliptic}, 
that is, $m$ has exactly one fixed point in $\d^-$, which lies in $\d$.
     \item $m$ is \emph{parabolic}, that is, $m$ has exactly one fixed point in $\d^-$,  which lies in $\t$.
         \item $m$ is \emph{hyperbolic}, that is, $m$ has exactly two fixed points in $\d^-$, which lie in $\t$.
 \end{enumerate}
We record the following lemma which gives a well-known criterion for the type of an automorphism in terms of the parameters $\alpha$ and $\tau$ in its representation as in equation \eqref{hyp.5}.
\begin{lem}\label{hyp.lem.5}
Let $m \in \aut \d\setminus \{\id{\d}\}$ and assume that $m =r_\tau \circ b_\alpha$ where $\alpha \in \d$  and $\tau \in \t$.
\begin{enumerate}[\rm (i)]
\item $m$ is elliptic$\iff |\tau - 1| > 2|\alpha|$.
\item $m$ is parabolic$\iff |\tau - 1| = 2|\alpha|$.
\item $m$ is hyperbolic$\iff |\tau-1| < 2|\alpha|$.
\end{enumerate}
\end{lem}
Note that this lemma implies that $b_{z_0}$ is hyperbolic for every $z_0 \in \d\setminus\{0\}$.
\begin{lem}\label{hyp.lem.7}
If $\alpha,\beta \in \d$ and $\beta \not=-\alpha$, then $b_\beta \circ b_\alpha$ is hyperbolic.
\end{lem}
\begin{proof}
Clearly $b_\beta \circ b_\alpha\neq \idd$.
By direct calculation,
\[
(b_\beta \circ b_\alpha) (z) = \frac{1+\bar\alpha\beta}{1+\bar\beta \alpha } \ \
\frac{z-\frac{\alpha +\beta}{1+\bar\alpha\beta}}
{1-\frac{\bar\alpha+\bar\beta}{1+\bar\beta\alpha}z}.
\]
Hence, Lemma \ref{hyp.lem.5} implies that $b_\beta \circ b_\alpha$ is hyperbolic if and only if
\begin{align*}
\Big|\frac{1+\bar\alpha\beta}{1+\bar\beta \alpha }-1\Big| <2\big|\frac{\alpha +\beta}{1+\bar\alpha\beta}\big|
\end{align*}
which simplifies to
\[
|\bar\alpha\beta -\bar\beta\alpha|<2|\alpha+\beta|.
\]
But, as $|\alpha|,|\beta|<1$ and $|\alpha+\beta|\neq 0$,
\begin{align*}
|\bar\alpha\beta -\bar\beta\alpha|&= |\bar\alpha(\alpha+\beta) -\beta(\bar\alpha +\bar\beta)|\\ 
&\le |\alpha||\alpha+\beta|+|\beta||\alpha+\beta|  <   2|\alpha+\beta|.
\end{align*}
\end{proof}
In the following definition we introduce a  class of automorphisms that plays a special role in the study of orthogonality.
\begin{defin}\label{hyp.def.2}
If $m \in \aut \d$ we say $m$ is \emph{irrotational} if in the representation of $m$ in equation \eqref{hyp.5}, $\tau=1$.
\end{defin}
\begin{remark}\label{irrothbol}
\rm
(1)  Lemma \ref{hyp.lem.5} implies that if $m$ is irrotational, then either $m=\id{\d}$ or $m$ is hyperbolic.  Lemma \ref{hyp.lem.7} now implies that if $m_1$ and $m_2$ are irrotational then either $m_2 =m_1^{-1}$ or $m _2\circ m_1$ is hyperbolic. \\
(2) Observe, by direct calculation, that for $m\in \aut \d$,
\[
m \text{ is irrotational } \iff m'(0) >0.
\]
\end{remark}
The following lemma characterizes when an automorphism is irrotational in terms of its fixed points.
\begin{lem}\label{hyp.lem.7a}
If $m \in \aut \d$, then $m$ is irrotational if and only if there exists $\eta \in \t$ such that $m(\eta) = \eta$ and $m(-\eta) =-\eta$.
\end{lem}
\begin{proof}
First assume that $m=b_\alpha$ is irrotational. When $\alpha=0$, \emph{every} point in $\t$ is a fixed point of $m$ and when $\alpha\not=0$, by direct computation the fixed points of $m$ are the roots of the equation $\eta^2=\alpha/\bar\alpha$. 

Conversely, assume that $\eta \in \t$, $m(\eta) = \eta$, and $m(-\eta) =-\eta$. By the representation for $m$ given in equation \eqref{hyp.5}, we have the equations
\begin{equation}\label{two-eq}
\tau\frac{\eta-\alpha}{1-\bar\alpha\eta}=\eta\ \ \text{ and }\ \
\tau\frac{-\eta-\alpha}{1+\bar\alpha\eta}=-\eta, \,
\end{equation}
which imply upon elimination of $\tau$ that $\eta^2 =\alpha/\bar\alpha$. By the first of equations \eqref{two-eq},
\[
\tau(\eta - \alpha) = \eta (1-\bar\alpha \eta)=\eta-\alpha,
\]
which implies that $\tau=1$.
\end{proof}
The following lemma examines a special point $\alpha \in \d$ that can be associated with an irrotational automorphism. This point $\alpha$ will play an essential role in our study of orthogonality in $G$.
\begin{lem}\label{hyp.lem.8}
Let $m \in \aut \d$. The following are equivalent. 
\begin{enumerate}[\rm (i)]
\item { $m$ is irrotational;}
\item {there exists $\alpha \in \d$ such that $m(\alpha)=-\alpha$ and $m'(\alpha) =1$;}
\item {there exists $\alpha \in \d$ such that $m=b_\alpha \circ b_\alpha$;}
\item { $-m$ is hyperbolic reflection about a point $\alpha$ in $\d$.}
\end{enumerate}
\end{lem}
\begin{proof}
(i) $\implies$ (ii). Let $m=b_\beta$ where $\beta\in\d$.  Lemma \ref{hyp.lem.5} implies that $-b_\beta$ is elliptic. Choose $\alpha \in \d$ such that $-b_\beta(\alpha) = \alpha$. We have
\begin{align*}
b_\beta'(\alpha)&= \frac{1-|\beta|^2}{(1-\bar\beta\alpha)^2}=\frac{(1-\bar\beta\alpha)+\bar\beta(\alpha-\beta)}{(1-\bar\beta\alpha)^2}\\ 
&=\frac{1}{1-\bar\beta\alpha} + \frac{\bar\beta}{1-\bar\beta\alpha}b_\beta(\alpha) =\frac{1}{1-\bar\beta\alpha} + \frac{\bar\beta}{1-\bar\beta\alpha}(-\alpha)=1.
\end{align*}
Hence (ii) holds.\\ 
(ii) $\implies$ (iii). One can show by Schur reduction that  $m = b_\alpha \circ r_\tau \circ b_\alpha$ for some $\tau \in \t$. But by the chain rule,
\[
1=m'(\alpha) =b_\alpha'(0)\cdot \tau \cdot b_\alpha'(\alpha)=(1-|\alpha|^2)\cdot \tau\cdot \frac{1}{1-|\alpha|^2}=\tau.
\]
(iii) $\implies$ (iv). By direct calculation, $-b_\alpha = b_{-\alpha} \circ r_{-1}$. Therefore, as $b_{-\alpha} = b_\alpha^{-1}$
\[
-b_\alpha \circ b_\alpha = b_\alpha^{-1}\circ r_{-1} \circ b_\alpha.
\]
Hence, as $b_\alpha^{-1}\circ r_{-1} \circ b_\alpha$ is hyperbolic reflection about $\alpha$, (iv) holds.\\ 
(iv) $\implies$ (i). If $-m$ is hyperbolic reflection about $\alpha$, then
\[
m = -b_\alpha^{-1}\circ r_{-1} \circ b_\alpha = b_\alpha \circ b_\alpha =b_\beta
\]
where
\[
\beta = \frac{2\alpha}{1+|\alpha|^2}.
\]
\end{proof}
Note that each of the
$\alpha$'s which appear in Conditions (ii), (iii), and (iv) of Lemma \ref{hyp.lem.8} are equal and are uniquely determined by $m$. Also note that the lemma implies that $b_\alpha \mapsto b_\alpha \circ b_\alpha$ is an injective map defined on the irrotational automorphisms onto the irrotational automorphisms. In particular, each irrotational automorphism has a unique irrotational square root with respect to composition.

\subsection{Purely balanced geodesics and royal points in $G$}
If $m \in \aut \d$, then the formula,
\beq\label{defhm}
h_m(z) =(z+m(z), zm(z)),\qquad z\in \d,
\eeq
defines a mapping $h_m:\d \to G$. \index{$ h_m$} It is easy to see that, when $m$ is elliptic,  $h_m$ is not injective, and thus cannot isometrically parametrize a geodesic. In all other cases, $h_m$ does parametrize a geodesic which we denote by $B_m$. The essentials of the correspondence between the properties of $m$ and $B_m$ are described in the next proposition. \index{$ B_m$} 
\begin{prop}\label{hyp.prop.5}
Let $m \in \aut \d$.
\begin{enumerate}[\rm (i)]
\item $m=\id{\d}$ $\iff$ $B_m=\royal$.
\item $m$ is hyperbolic $\iff$ $B_m$ is purely balanced.
\item $m$ is parabolic $\iff$ $B_m$ is exceptional.
\end{enumerate}
\end{prop}
\begin{proof}
This proposition  is proved in \cite[Theorem 7.8]{aly2016}.
\end{proof}

\begin{lem}\label{mqinv}
If $m,q$ are hyperbolic automorphisms of $\d$ and $B_m=B_q$ then either $q=m$ or $q=m\inv$.
\end{lem}
\begin{proof}
Let the fixed points of $m, q$ be $\eta=\{\eta_1,\eta_2\}$ and $\zeta=\{\zeta_1,\zeta_2\}$ respectively.
Then $h_q\inv\circ h_m$ is an automorphism, $\ups$ say, of $\d$.  Thus $h_m=h_q\circ \ups$, that is, for all $z\in\d$,
\[
(z+m(z),zm(z))= h_q(\ups(z))=(\ups(z)+q\circ \ups(z),\ups(z) q\circ\ups(z)).
\]
Therefore, for all $z$, either
\[
\ups(z)=z \quad \mbox{ and } \quad  q\circ \ups(z) = m(z)
\]
or
\[
q\circ \ups(z) = z \quad  \mbox{ and }\quad  \ups(z) = m(z).
\]
Thus either $\ups=\idd $ and $q=m$, or $\ups=m$ and $q=m\inv$.
\end{proof}

For hyperbolic $m \in \aut \d$ there are  important relationships between the fixed points of $m$ and properties of the corresponding purely balanced geodesic $B_m$. We describe three of these in Lemma \ref{hyp.lem.20}, Corollary \ref{hyp.lem.30} and Proposition \ref{dis.prop.40} below.

Recall that a geodesic $D$ was defined to be purely balanced if $D^- \cap \royal^-$ consists of exactly two points, both of which lie in ${\mathrm E}$.
\begin{lem}\label{hyp.lem.20}
If $m$ is hyperbolic and $\eta_1$ and $\eta_2$ are the distinct fixed points of $m$ in $\t$, then
\[
B_m^- \cap \royal^- =\big\{(2\eta_1\eta_1^2),(2\eta_2,\eta_2^2)\big\}.
\]
\end{lem}
\begin{proof}
If $m(\eta) =\eta$, then $h_m(\eta) =(2\eta,\eta^2)$. Therefore,
\[
(2\eta_1\eta_1^2),(2\eta_2,\eta_2^2)\in  B_m^- \cap \royal^-.
\]

Conversely, if $\lambda \in B_m^- \cap \royal^-$, then since $\lambda \in B_m^-$, there exists $\zeta \in \t$ such that
\[
\lambda = (\zeta +m(\zeta), \zeta m(\zeta)),
\]
and since $\lambda \in \royal^-$, there exists $\eta \in \t$ such that
\[
\lambda =(2\eta, \eta^2).
\]
These equations imply that
\[
\zeta +m(\zeta) = 2\eta\ \ \text{ and }\zeta m(\zeta)=\eta^2.\ \
\]
Hence, for all $x\in\c$,
\begin{align*}
(x-\zeta)(x-m(\zeta))&= x^2 -(\zeta +m(\zeta))x + \zeta m(\zeta)\\
&=x^2 -2\eta x +\eta^2 =(x-\eta)^2,
\end{align*}
which implies that $\zeta=\eta$ and so $m(\eta) =\eta$.
\end{proof}
A second relationship between $B_m$ and the fixed points of $m$ involves the qualitative nature of the solutions to the Carath\'eodory extremal problem.  We say, for a domain $\Omega$, that a set $\mathcal C$ of holomorphic maps from $\Omega$ to $\d$ is a {\em universal set for the Carath\'eodory extremal problem on $\Omega$} if, for every $\la\in \Omega$ and every nonzero vector $v\in T_\la \Omega$, there exists $\Phi \in \mathcal C$ such that the supremum in equation \eqref{defCmetric}, when $U = \Omega$, is attained at $F=\Phi$.

For $\omega \in \t$, define a  holomorphic function $\Phi_\omega$ on $G$ by the formula,
\[
\Phi_\omega(s,p) =\frac{2\omega p -s}{2-\omega s},\qquad \fa (s,p) \in G.
\]

The following properties of $\Phi_\omega, \ \omega\in\t,$ are established in \cite[Theorem 2.1 and Corollary 4.3]{ay2004}.
\begin{prop}\label{hyp.prop.10}
For every $\omega\in \t$,  $\Phi_\omega$ is a holomorphic map from $G$ to $\d$. Furthermore, the set 
$\set{\Phi_\omega}{\omega \in \t}$ is universal for the Carath\'eodory extremal problem on $G$.
\end{prop}

The following result gives a criterion for a geodesic in $G$ to be purely balanced in terms of a qualitative property of solutions of the Carath\'eodory extremal problem.  It is contained in 
\cite[Theorem 7.1(iii)]{aly2016}.
\begin{prop}
Let $D$ be a geodesic in $G$, let $\la\in D$ and let $v$ be a nonzero vector in $T_\la D$.  The geodesic $D$ is purely balanced if and only if there are two distinct points $\omega_1,\omega_2 \in \t$ such that the supremum in equation \eqref{defCmetric}, when $U = G$,  is attained at $F=\Phi_\omega$ precisely when either $\omega=\omega_1$ or $\omega=\omega_2$.  Furthermore, if $D$ is purely balanced and $\omega_1,\omega_2$ are as in the preceding sentence, then the royal points of $D$ are $(2\bar\omega_1,\bar\omega_1^2)$ and $(2\bar\omega_2,\bar\omega_2^2)$.
\end{prop}

The second promised relationship between a purely balanced geodesic and its royal points is contained in the following statement.
\begin{corollary}\label{hyp.lem.30}
Let $m$ be a hyperbolic automorphism of $\D$, let $\eta_1$ and $\eta_2$ be the distinct fixed points of $m$ in $\t$, let $\la\in B_m$ and let $v\in T_\la B_m$. Then $\Phi_\omega$ is a solution of the Carath\'eodory extremal problem on $G$ at the point $\la$ in direction $v$ if and only if either
$\omega=\bar\eta_1$ or $\omega=\bar\eta_2$.
\end{corollary}

\subsection{Purely balanced geodesics in standard position in $G$}
In this subsection we set out a canonical form for purely balanced geodesics. We accomplish this goal by cleanly parametrizing the purely balanced geodesics that pass through a point of the form $(0,p)$ where $p<0$.
\begin{lem}\label{hyp.lem.40}
Let $B$ be a purely balanced (or exceptional) geodesic and let $\sigma\in (0,1)$.
\[
(0,-\sigma^2) \in B
\]
 if and only if there exists $m \in \aut \d$ satisfying
\begin{enumerate}[\rm (i)]
\item $m$ is hyperbolic (or parabolic respectively),
\item $m(\sigma) = -\sigma$, and
\item $B=B_m$.
\end{enumerate}
\end{lem}
\begin{proof}
If conditions (i), (ii) and (iii) hold, then
\[
(0,-\sigma^2) = \big(\sigma +(-\sigma),\sigma(-\sigma)\big)=
h_m(\sigma) \in B_m=B.
\]
Furthermore, Proposition \ref{hyp.prop.5} implies that $B$ is purely balanced (or exceptional respectively).

Conversely, assume that $B$ is purely balanced (or exceptional) and $(0,-\sigma^2) \in B$. Using Proposition \ref{hyp.prop.5} we may choose a hyperbolic (or exceptional respectively) $b\in \aut \d$ such that $B=B_b$.
Choose $z_0 \in \d$ such that $h_b(z_0)=(0,-\sigma^2)$. Then
\[
b(z_0) = -z_0\ \ \text{ and }\ \ z_0^2 =\sigma^2.
\]
If $z_0=\sigma$, then the conditions (i), (ii), and (iii) follow if we set $m=b$. Otherwise, if $z_0 =-\sigma$, the conclusions follow if we set $m=b^{-1}$.
\end{proof}
\begin{lem}\label{hyp.lem.50}
If $\sigma >0$ and $m \in \aut \d$, then $m(\sigma) = -\sigma$ if and only if $m=b_\sigma \circ r_\tau \circ b_\sigma$ for some $\tau \in \t$. Furthermore, if $m(\sigma) =-\sigma$, then $m$ is hyperbolic if and only if
\be\label{hyp.30}
\left|\frac{\tau-1}{\tau +1}\right| < \frac {2\sigma}{1-\sigma^2},
\ee
and $m$ is parabolic if and only if
\be\label{hyp.40}
\left|\frac{\tau-1}{\tau +1}\right| =\frac {2\sigma}{1-\sigma^2}.
\ee
\end{lem}
\begin{proof}

\begin{align*}
m(\sigma) =-\sigma &\iff  m (b_\sigma^{-1}(0)) =b_\sigma(0) \\
&\iff\big(b_{\sigma}^{-1} \circ m \circ b_\sigma^{-1}\big)(0) = 0\\
&\iff  b_{\sigma}^{-1} \circ m \circ b_\sigma^{-1}=r_\tau\quad \fs \tau\in \t\\
&\iff m=b_\sigma \circ r_\tau \circ b_\sigma\quad \fs \tau\in\t,
\end{align*}
which proves the first assertion of the lemma.

To prove the second assertion, by direct calculation,
\[
m(z) = \frac{\tau +\sigma^2}{1+\tau \sigma^2}\ \
\frac{z-\frac{\sigma(\tau+1)}{\tau+\sigma^2}}
{1- \frac{\sigma(\tau+1)}{1+\tau \sigma^2} z}.
\]
Consequently, by Lemma \ref{hyp.lem.5} (iii), $m$ is hyperbolic if and only if
\begin{align*}
& \left|\frac{\tau + \sigma^2}{1+\tau \sigma^2} - 1\right| <2\ \left|\frac{\sigma(\tau+1)}{1+\tau \sigma^2}\right|\\ 
\iff\ \  & \left|\tau +\sigma^2 -(1 +\tau \sigma^2)\right| <2\  \left|\sigma(\tau+1)\right|\\ 
\iff\ \ & \left|\frac{\tau-1}{\tau +1}\right| < \frac {2\sigma}{1-\sigma^2},
\end{align*}
that is,  inequality \eqref{hyp.30} holds.

 A similar calculation using Lemma \ref{hyp.lem.5} (ii) shows that $m$ is parabolic if and only if equation \eqref{hyp.40} holds.
\end{proof}
\subsection{The curve of purely balanced directions}
If $\tau_1$ and $\tau_2$ are in $\t$, we denote by $(\tau_1,\tau_2)$ the open segment of points in $\t$ extending from  $\tau_1$ to $\tau_2$ in the counterclockwise direction. In similar fashion we may define half open and closed segments in the circle.  For $\sigma \in (0,1)$,  let $\tau_\sigma^+$ denote the unique point in $\t$ that satisfies equation  \eqref{hyp.40} and $\im \tau_\sigma^+ >0$ and let $\tau_\sigma^-$ denote the unique point in $\t$ satisfying equation \eqref{hyp.40} and $\im \tau_\sigma^- <0$.
When $\sigma \in (0,1)$ and $\tau \in \t$, we define $m_{\sigma,\tau}$ by 
\be\label{defmst}
m_{\sigma,\tau}=b_\sigma \circ r_\tau \circ b_\sigma.
\ee
\index{$m_{\si,\tau}$}

\begin{lem}\label{parampb}
The map $(\si,\tau) \mapsto B_{m_{\si,\tau}}$ is injective from  the set
\[
X\df \{(\si,\tau): 0 < \si < 1, \, \tau \in (\tau_\si^-, \tau_\si^+)\}
\]
into the set of purely balanced geodesics in $G$.
\end{lem}
\begin{proof}
It follows from Lemmas \ref{hyp.lem.50} and \ref{hyp.lem.40} that, for $(\si,\tau)\in X$, the geodesic $B_{m_{\si,\tau}}$ is purely balanced.

Suppose that two points $(\si,\tau)$ and $(\si',\tau')$  in $X$ give rise to the same  purely balanced geodesic $B$.  Notice firstly that
$B$ only meets the flat geodesic $F^0= \{0\}\times \d$ in a single point, else, by the uniqueness of geodesics through any pair of distinct points, we have $B=F_0$, contrary to the hypothesis that $B$ is purely balanced.  By Lemma \ref{hyp.lem.40}, both $(0,-\si^2)$ and $(0,-(\si')^2)$ lie in $B$, hence coincide.  Since $\si,\si'$ are positive, it follows that $\si=\si'$.

To see that $\tau=\tau'$, calculate a tangent vector to $B$ at $(0,-\si^2)$.  We find that
\[
h_{m_{\si,\tau}}'(\si) = (1+\tau, -\si(1-\tau)).
\]
It follows that, for some nonzero complex number $c$,
\[
(1+\tau, -\si(1-\tau)) = c(1+\tau', -\si(1-\tau')),
\]
from which it follows that $\tau=\tau'$.
\end{proof}

With the notations of the previous paragraph, the following proposition is an immediate consequence of
Lemmas \ref{hyp.lem.40} and \ref{hyp.lem.50}.
\begin{prop}\label{hyp.prop.20}
Let $\sigma \in (0,1)$ and let $D$ be a geodesic in $G$.  Then $D$ is a purely balanced geodesic passing through $(0,-\sigma^2)$ if and only if there exists $\tau \in (\tau_\sigma^-,\tau_\sigma^+)$ such that $D=B_{m_{\sigma,\tau}}$. Moreover $D$ is an exceptional geodesic passing through $(0,-\sigma^2)$ if and only if $D=B_{m_{\sigma,\tau}}$ where either $\tau=\tau_\sigma^+$ or $\tau=\tau_\sigma^-$.
\end{prop}
We refer to geodesics of the form $B_{m_{\sigma,\tau}}$ as being in \emph{standard position}. The following proposition states that any purely balanced or exceptional geodesic with a point prescribed in it can be moved to a geodesic in standard position by an automorphism of $G$.
\begin{prop}\label{hyp.prop.30}
Let $D$ be a geodesic in $G$ that is purely balanced (or exceptional) and let $\lambda \in D$.  There exist $\sigma \in (0,1)$,
$\tau\in (\tau_\sigma^-,\tau_\sigma^+)$
(or $\tau \in \{\tau_\sigma^-,\tau_\sigma^+\}$ respectively),
 and $\gamma \in \aut G$ such that $\gamma(\lambda) =(0,-\sigma^2)$ and $\gamma(D) = B_{m_{\sigma,\tau}}$.
\end{prop}
\begin{proof}
Let $\lambda = \pi(z)$ where $z=(z_1,z_2) \in \d^2$. Since $\lambda$ is in the purely balanced or exceptional  geodesic $D$, $\lambda \not\in \royal$. Consequently, $z_1\not=z_2$ and by the intermediate value theorem we may choose $\sigma \in (0,1)$ such that
\[
d(\sigma, -\sigma) = d (z_1,z_2).
\]
Choose $b \in \aut \d$ such that
\[
b(z_1) = \sigma\ \ \text{ and }\ \ b(z_2) = -\sigma.
\]
By Proposition \ref{pre.prop.5},  $\gamma_b \in \aut G$. Furthermore, by formula \eqref{pre.5},
\[
\gamma_b(\lambda)=(b(z_1)+b(z_2),b(z_1)b(z_2))=(0,-\sigma^2).
\]
As $B$ is assumed purely balanced or exceptional, so also is $\gamma_b(B)$, and $(0,-\si^2)\in \ga_b(\d)$. Therefore, the proposition follows by  Proposition \ref{hyp.prop.20}.
\end{proof}
\begin{corollary}\label{hyp.cor.10}
 If $\lambda \in \royal$, then there are no purely balanced or exceptional directions at $\lambda$. If $\lambda \in G \setminus\royal$, then there are exactly two exceptional directions at $\la$, and the purely balanced directions at $\la$ form a simple smooth curve in $\cptwo$ connecting the exceptional directions.
\end{corollary}
\begin{proof}
As no purely balanced or exceptional geodesic meets $\royal$, there are no purely balanced or exceptional directions at points in $\royal$. To prove the second assertion of the corollary, by Proposition \ref{hyp.prop.30}, it suffices to prove the case when $\lambda=(0,-\sigma^2)$ for some $\sigma \in (0,1)$. But this case is an immediate consequence of Proposition \ref{hyp.prop.20}.
\end{proof}

%% file: orthogonality.tex

 \section{Orthogonality in $G$}\label{orthogonality}
In this section we shall study the geodesics that are orthogonal to a fixed flat geodesic.
\subsection{Flat geodesics}\label{flgeo}
Recall that the flat geodesics of $G$ were defined in Subsection \ref{5types} to be the geodesics $D$ that have a unique royal point $\la$ and furthermore are such that $\la\in\calr$.
It can be shown that flat geodesics are truly `flat': they have
 the form\footnote{In earlier papers we {\em defined} the flat geodesics to be the $F^\beta$ and showed that they are characterized by the property that they have a single royal point, which lies in $\calr$.  Here, we reverse the process in order to bring out the geometric nature of flatness.}
\be\label{flatbeta}
F^\beta \df \{(\beta+\bar\beta z,z): z\in\d\}
\ee
 for some $\beta \in\d$.  
\index{$F^\beta$}
\index{geodesic!flat}
Let us check that, for every $\beta\in\d$, this set is indeed a flat geodesic of $G$ according  to the definition.  Firstly, if $(s,p) = (\beta+\bar\beta z,z)$ for some $z, \beta \in \d$, then $z=p$ and
\[
s-\bar s p= \beta+\bar\beta p-(\bar\beta+\beta \bar p) = \beta(1-|p|^2),
\]
whence
\[
|s-\bar s p| < 1-|p|^2.
\]
Hence, by the criterion \eqref{schurcrit}, $(s,p)\in G$.  Thus $F^\beta \subset G$.  Secondly, $F^\beta$ is the range of the analytic disc $h(z)=(\beta+\bar\beta z,z)$ in $G$, and $h$ has an obvious holomorphic left inverse, to wit, the second co-ordinate function.  Thus $F^\beta$ is a geodesic of $G$.  

It remains to show that $F^\beta$ meets $\calr^-$ at exactly one point $\la$, and that $\la\in\calr$.   Indeed, $(s,p)=h(z)\in\calr^-$ if and only if
\[
0=s^2-4p=(\beta+\bar\beta z)^2 -4z = \bar\beta^2 z^2+2(|\beta|^2-2)z + \beta^2.
\]
If $\beta=0$ then this equation for $z$ has the unique root $z=0$, so that $F^0$ meets $\calr^-$ at the unique point $(0,0)$.  Otherwise,
the quadratic equation has two roots, whose product has modulus one and whose sum has modulus greater than two.  It follows that the equation has exactly one root in $\d$ and none in $\t$.  Thus, in either case, $F^\beta$ has a unique royal point, and that point lies in $\calr$, which is to say that $F^\beta$ is a flat geodesic.
\begin{prop}\label{permuteflat}
$\aut G$ permutes the set of flat geodesics of $G$.   Indeed, for any $\beta\in\d$,
\[
\ga_{b_\al}(F^\beta)=F^0,
\]
where 
\be\label{hmidpt}
\al=\frac{\beta}{1+\sqrt{1-|\beta|^2}}
\ee
is the hyperbolic mid-point of $0$ and $\beta$ in $\d$.
\end{prop}
\begin{proof}
By Corollary \ref{extends}(3), every automorphism of $G$ maps every flat geodesic to a flat geodesic.

  $\aut G$ acts transitively on $\calr$, by reason of the equation
\[
\ga_m(2z,z^2) = (2m(z),m(z)^2)
\]
for every $m\in\aut \d$ and $z\in\d$.  Hence, if $\la, \mu \in \calr,$ then we can find $m\in \aut\d$ such that $\ga_m(\la)=\mu$.  It follows that $\ga_m(F_\la)$ is a flat geodesic that contains $\mu\in\calr$.  Since there is a {\em unique} flat geodesic through any point of $G$, we infer that  $\ga_m(F_\la)=F_\mu$.  Thus $\aut G$ permutes the flat geodesics of $G$.
 
It is not too arduous to calculate that, if $(s,p)\in F^\beta$ and $\al$ is given by equation \eqref{hmidpt}, then
\[
\ga_{b_\al}(s,p)= \left( 0, \frac{p-\al s + \al^2}{1-\bar\al s + \bar\al^2 p}\right)
\]
which is in $F^0$.
\end{proof}

\subsection{The sharp direction}\label{sharpdir}
Since every orbit of $\aut G$ meets the flat geodesic $F^0$ and the sharp direction is covariant,
the sharp direction at any point in $G$ can be derived from the following statement.

\begin{prop}\label{sharpF0}
For any
$\la\in F^0$,
\be\label{sharpspecial}
\sharp_\la=\bpm 1 \\ 0 \epm \c.
\ee
\end{prop}
\begin{proof}
Consider the point $\la=(0,p)\in F^0$.  For any nonzero $\al\in\c$, let
\[
f_\al(t)=\ga_{b_{t\al}}(\la) \qquad \mbox{ for } t\in\r \mbox{ such that }|t\al| < 1.
\]
Then $f_\al(t)$ describes a smooth path in $\orb(\la)$ such that $f_\al(0)=\la$, and hence $f_\al'(0) \in T_\la\orb(\la)$.

Let $w\in\d$ satisfy $w^2=-p$, so that $\la=\pi(-w,w)$.  A short calculation shows that
\begin{align*}
\ga_{b_\al}(\la) &= \pi(b_\al(w),b_\al(-w)) \\
	&=\frac{1}{1+\bar\al^2p} \left(-2(\al+\bar\al p), p+\al^2\right).
\end{align*}
Replace $\al$ by $t\al$, differentiate with respect to $t$ and set $t=0$ to obtain
\[
f_\al'(0)=\left( -2(\al+\bar\al p), 0 \right) \quad \in T_\la\orb(\la)
\]
for every $\al \in \c\setminus\{0\}$.  By suitable choices of $\al$ we deduce that both $(1,0)$ and $(i,0)$ lie in $T_\la\orb(\la)$.  Thus the complex linear span of $(1,0)$ is contained in $T_\la\orb(\la)$, and statement \eqref{sharpspecial} follows.
\end{proof}

\begin{corollary}
For every $\la\in G, \; \sharp_\la \neq \flat_\la$.
\end{corollary}
\begin{proof}
It is enough to prove it for every $\la\in F^0$.  For such $\la$,
\[
\sharp_\la= \bpm 1 \\ 0 \epm\c, \quad \flat_\la=\bpm 0\\ 1 \epm\c.
\]
\end{proof}
It follows that the tangent bundle $TG$ is the direct sum of the sharp bundle $\la\mapsto \sharp_\la$ and the flat bundle $\la\mapsto \flat_\la$.

We remark that a formula for $\sharp_\la$ in general is given in \cite[Proposition 1.16]{aly2017}.
If $\lambda=(s,p) \in F^\beta$, then 
\[
\sharp_\lambda = \begin{pmatrix}1 \\  \ds  \frac{\beta-\half s}{1-\half\bar\beta s} \end{pmatrix}\c.
\]

\subsection{Orthogonality and irrotational automorphisms}
Recall that a geodesic $D$ is said to be orthogonal to a flat geodesic $F$ if $D$ meets $F$ at a point $\lambda$ and $T_\lambda D = \sharp_\lambda$.
In view of Proposition \ref{permuteflat}, in studying any flat geodesic $F$, one may often reduce to the case that 
\[
F=F^0= \{(0,p): p\in\d\},
\]
the flat geodesic passing through the origin. 

\begin{prop}\label{sharp.prop.10}
A geodesic $D$ is orthogonal to $F^0$ if and only if $D=B_m$ for some irrotational $m \in \aut \d$.
\end{prop}
\begin{proof}
First assume that $D \perp F^0$. Choose $\lambda \in D\cap F^0$. Let $\lambda=(0,p_0)$ and choose $\alpha \in \d$ such that $\alpha^2 =-p_0$. Let $m\in \aut \d$ where $-m$ is hyperbolic reflection about $\alpha$. Observe that Lemma \ref{hyp.lem.8} implies that $m$ is irrotational.

To see that $D=B_m$ note that Condition (ii) in Lemma \ref{hyp.lem.8} implies that $m(\alpha)=-\alpha$, which implies that $h_m(\alpha) =\lambda$. Therefore, $\lambda \in B_m$. Also, Condition (ii)  in Lemma \ref{hyp.lem.8} implies that
$m'(\alpha) =1$. Therefore,
\begin{align*}
h_m'(\alpha)=(1 +m'(\alpha),m(\alpha)+\alpha m'(\alpha))=(1+1,-\alpha+\alpha) = (2,0),
\end{align*}
which implies that $T_\la B_m=\sharp_\lambda=T_\lambda D$. As $D$ and $B_m$ both pass through the point $\lambda$ and have the same tangent space at $\lambda$, it follows from Theorem \ref{pre.thm.10} that $D=B_m$.

Now assume that $m$ is irrotational. We need to show that $B_m \perp F^0$.  Condition (ii) in Lemma \ref{hyp.lem.8} implies that there exists $\alpha \in \d$ such that $m(\alpha)=-\alpha$. As,
\[
h_m(\alpha)=(0,-\alpha^2) \in F^0,
\]
it follows that $B_m$ meets $F^0$ at the point $\lambda=(0,-\alpha^2)$.
In addition,  Condition (ii) in Lemma \ref{hyp.lem.8} implies that $m'(\alpha)=1$. Therefore, calculating as before, we have
\[
h_m'(\alpha) =(2,0),
\]
which implies that $T_\la B_m = \sharp_\lambda$. As  $B_m$ meets $F^0$ at the point $\lambda$ and $T_\la B_m = \sharp_\lambda$, it follows that $B_m$ is orthogonal to $F^0$.
\end{proof}
\begin{corollary}\label{sharp.cor.20}
If $\mu \in \royal$, then $\sharp_\mu=T_\mu\royal$. If $\mu \in G \setminus \royal$, then the geodesic through $\mu$ with direction $\sharp_\mu$ is purely balanced.
\end{corollary}
\begin{proof}
Fix $\mu \in G$, let $F$ be the flat geodesic passing through $\mu$, and let $D$ be the geodesic such that $\mu \in D$ and $T_\mu D =\sharp_\mu$. We wish to show that $D=\royal$ when $\mu \in \royal$, and otherwise, that $D$ is purely balanced.

Choose $\gamma \in \aut G$ such that $\lambda =\gamma(\mu)\in F^0$. As $\gamma$ is an automorphism and $D \perp F$, $\gamma(D) \perp F^0$. Therefore, by Proposition \ref{sharp.prop.10}, there exists an irrotational $m$ such that $\gamma(D) = B_m$.

As $m$ is irrotational, by Proposition \ref{hyp.prop.5} there are only two possibilities:
\begin{enumerate}[(i)]
\item $m=\id{\d}$ and $B_m=\royal$, and
\item $m$ is hyperbolic and $B_m$ is purely balanced.
\end{enumerate}
In case (i), since $\gamma$ fixes $\royal$, $\mu \in \royal$ and $\sharp_\mu = T_\mu\royal$. In case (ii), $\mu \not\in \royal$, and $\sharp_\mu$ is purely balanced.
\end{proof}
The following result is essentially Corollary \ref{sharp.cor.20} stated in different language.
\begin{corollary}\label{sharp.cor.30}
The royal geodesic is orthogonal to every flat geodesic. If a geodesic $D$ is orthogonal to a flat geodesic, then either $D=\royal$ or $D$ is purely balanced.
\end{corollary}
The following result implies that at points not on the royal geodesic, the sharp direction is the midpoint of the curve of purely balanced directions.
\begin{corollary}\label{sharp.cor.40}
For the curve of purely balanced geodesics $B_\tau = B_{m_{\sigma,\tau}}$ constructed in Proposition {\rm \ref{hyp.prop.20}}, $T_{(0,-\si^2)}B_\tau = \sharp_{(0,-\sigma^2)}$ if and only if $\tau=1$.
\end{corollary}
\begin{proof}
Proposition \ref{sharp.prop.10} implies that $T_{(0,-\sigma^2)}  B_\tau = \sharp_{(0,-\sigma^2)}$ if and only if $m_{\sigma,\tau}$ is irrotational. But Lemma \ref{hyp.lem.8} implies that $m_{\sigma,\tau}$ is irrotational if and only if $\tau=m_{\sigma,\tau}'(\sigma)=1$.
\end{proof}

\subsection{Foliation of $G$ by orthogonal geodesics}\label{fol_pure_bal}
If $F$ is a geodesic in $\r^2$ (that is, a line), then $\r^2$ is foliated by the geodesics orthogonal to $F$. Theorem \ref{sharp.thm.10} below is an analogous result in $G$. 

\begin{lem}\label{sharp.lem.90}
If $F$ is a flat geodesic in $G$, and $D_1$ and $D_2$ are geodesics that are orthogonal to $F$, then either $D_1=D_2$ or $D_1 \cap D_2=\varnothing$.
\end{lem}
\begin{proof}
As $\aut G$ acts transitively on flat geodesics,
 we may assume without loss of generality that $F=F^0$.
Let $D_1$ and $D_2$ be geodesics that are orthogonal to $F$. By Proposition \ref{sharp.prop.10}, there exist irrotational automorphisms $m_1$ and $m_2$ of $\d$ such that $D_1 = B_{m_1}$ and $D_2=B_{m_2}$.

If $B_{m_1}\cap B_{m_2} \not=\varnothing$, then there exist $z_1$ and $z_2$ in $\d$ such that $h_{m_1}(z_1) =h_{m_2}(z_2)$,
equations which imply either
\be\label{sharp.80}
z_1= z_2\ \ \text{ and }\ \  m_1(z_1) = m_2(z_2)
\ee
or
\be\label{sharp.90}
z_1= m_2(z_2)\ \ \text{ and }\ \  m_1(z_1) = z_2.
\ee

If equations \eqref{sharp.80} hold then $m_1(z_1) = m_2(z_1)$, or equivalently, $(m_2^{-1}\circ m_1) (z_1) = z_1$. Hence, $m_2^{-1} \circ m_1$ is elliptic, and Remark \ref{irrothbol}(1) implies that $m_1=m_2$. Thus $D_1=D_2$.

If equations \eqref{sharp.90} hold, then $(m_2\circ m_1)(z_1)=z_1$ and we see that $m_2\circ m_1$ is elliptic. Then Remark \ref{irrothbol}(1) implies that $m_1=m_2^{-1}$, and again,  $D_1=D_2$.
\end{proof}
\begin{lem}\label{sharp.lem.100}
Fix a flat geodesic $F$ and a point $\mu \in G\setminus F$. There exists a purely balanced geodesic $D$ such that $D$ is orthogonal to $F$ and $\mu \in D$.
\end{lem}
\begin{proof}
Without loss of generality we may assume that $F=F^0$. Choose $z_1,z_2 \in \d$ such that
\[
\mu =\pi(z_1,z_2)=(z_1+z_2,z_1z_2).
\]
Represent $\pi(z_1,-z_2)$ in flat coordinates,
\[
\pi(z_1,-z_2)=(z_1-z_2,-z_1z_2) =(\beta +\bar\beta p_0,p_0),
\]
where $\beta,p_0 \in \d$.
We have
\[
z_1-z_2=\beta -\bar\beta z_1z_2,
\]
which implies that
\[
z_2=\frac{z_1-\beta}{1-\bar\beta z_1}=b_\beta(z_1).
\]
If we set $D=B_{b_\beta}$, then, as $h_{b_\beta}(z_1)=\mu$, we can assert that $\mu \in D$. Furthermore, as $b_\beta$ is irrotational, Proposition \ref{sharp.prop.10} guarantees that $D$ is orthogonal to $F$.
\end{proof}
By combining the previous two lemmas we obtain the following theorem.
\begin{thm}\label{sharp.thm.10}
If $F$ is a flat geodesic, then $G$ is foliated by the geodesics in $G$ that are orthogonal to $F$.
\end{thm}

%% file: distinguished.tex

\section{Distinguished geodesics in $G$}\label{distinguished}

In the previous section we studied the purely balanced geodesics that are orthogonal to a fixed flat geodesic. In this section, in contrast, we
fix a purely balanced geodesic $D$ and study  the flat geodesics $F$ such that $D$ is  orthogonal to $F$. This leads to the remarkable discovery that such flat geodesics are naturally parametrized by a real geodesic in $D$.
\begin{defin}\label{dis.def.10}
For any geodesic $D$ we say that a point $\la\in D$ is a {\em sharp point in $D$} if $D$ is orthogonal to the flat geodesic passing through $\lambda$.   The set of sharp points in $D$ is denoted by $\sharp(D)$.
\end{defin}
\index{sharp!point}
\index{$\sharp(D)$}
Thus $\sharp(D)=\{\la\in D: T_\la D=\sharp_\la\}$.
 Corollary \ref{sharp.cor.20} implies that $\sharp(\royal) =\royal$ and that, if $D\neq \calr$, then $\sharp(D)$ is non-empty only if $D$ is purely balanced.  We therefore restrict attention to the case that $D$ is a purely balanced geodesic $B_m$.

Proposition \ref{dis.prop.40} below, which gives a description of the sharp points in a purely balanced geodesic $B$,  represents the third promised relationship between $B= B_m$ and the fixed points of $m$.

If $B$ is a purely balanced geodesic, we may define a real geodesic in $B$ in the following manner. By  Proposition \ref{hyp.prop.5}, there exists hyperbolic $m\in\aut\d$ such that $B=B_m$.  Let $\eta=\{\eta_1,\eta_2\}$ denote the set of fixed points of  $m$ and let $C_\eta$ denote the  real hyperbolic geodesic in $\d$ that has $\eta_1$ and $\eta_2$ as endpoints. Finally, we define $\Xi_B \subset B$ by
\[
\Xi_B = h_m(C_\eta),
\]
where, for $ z\in \d$, $h_m(z) = (z + m(z), z m(z))$.
\begin{lem}\label{dis.lem.10} 
If $B$ is a purely balanced geodesic and $B=B_m$, then the definition of $\Xi_B$ does not depend on the choice of $m$.
\end{lem}
\begin{proof}
Suppose that $m,q \in\aut\d$ are such that $B_m=B=B_q$.
By Lemma \ref{mqinv},  either $q=m$ or  $q=m\inv$.
In the former case it is immediate that $m$ and $q$ yield the same curve $\Xi_B$.
In the latter case,  since the fixed points $\zeta$ of $m\inv$ coincide with those $\eta$ of $m$, we have $\eta=\zeta$,
and so $C_\eta=C_\zeta$.   Moreover, since $m$ is a hyperbolic isometry of $\d$, $m$ maps geodesics to geodesics in $\d$, and so $m(C_\eta)=C_\eta=C_\zeta$.  For any $z\in\d$, we have $h_{m\inv}(m(z))= h_{m}(z)$, and therefore
\[
h_q(C_\zeta)=h_{m\inv}(m(C_\eta))=h_m(C_\eta) = \Xi_B.
\]
\end{proof}
\begin{prop}\label{dis.prop.40}
If $B$ is a purely balanced geodesic, then $\sharp(B) = \Xi_B$.
\end{prop}
\begin{proof}
We first show that $\sharp(B)  \subseteq \Xi_B$. Fix $\lambda \in \sharp(B)$. By Proposition \ref{hyp.prop.30}, we may assume that $B$ is in standard position, $\lambda=(0,-\sigma^2)$ and $B=B_{m_{\sigma,\tau}}$ for some $\sigma \in (0,1)$ and $\tau \in \t$. Furthermore, as $\lambda \in \sharp(B)$, Corollary \ref{sharp.cor.40} implies that $\tau =1$. Since the fixed points of $m_{\sigma,1}$ are $\pm 1$, $C_\eta =(-1,1)$. Therefore, as $\sigma \in (-1,1)$,
\[
\lambda = (0,-\sigma^2) = h_{m_{\sigma,1}}(\sigma) \in h_{m_{\sigma,1}}(C_\eta) =\Xi_B.
\]

We now turn to the proof that $\Xi_B \subseteq \sharp(B)$. Fix a purely balanced geodesic $B$ with royal points $(2\xi,\xi^2)$ and $(2\eta,\eta^2)$. Assume that $B$ is parametrized as in equation \eqref{defhm}, so that $B=B_m$ for some hyperbolic $m\in \aut \d$, with fixed points at $\xi$ and $\eta$ in $\t$.

\begin{center}
\includegraphics{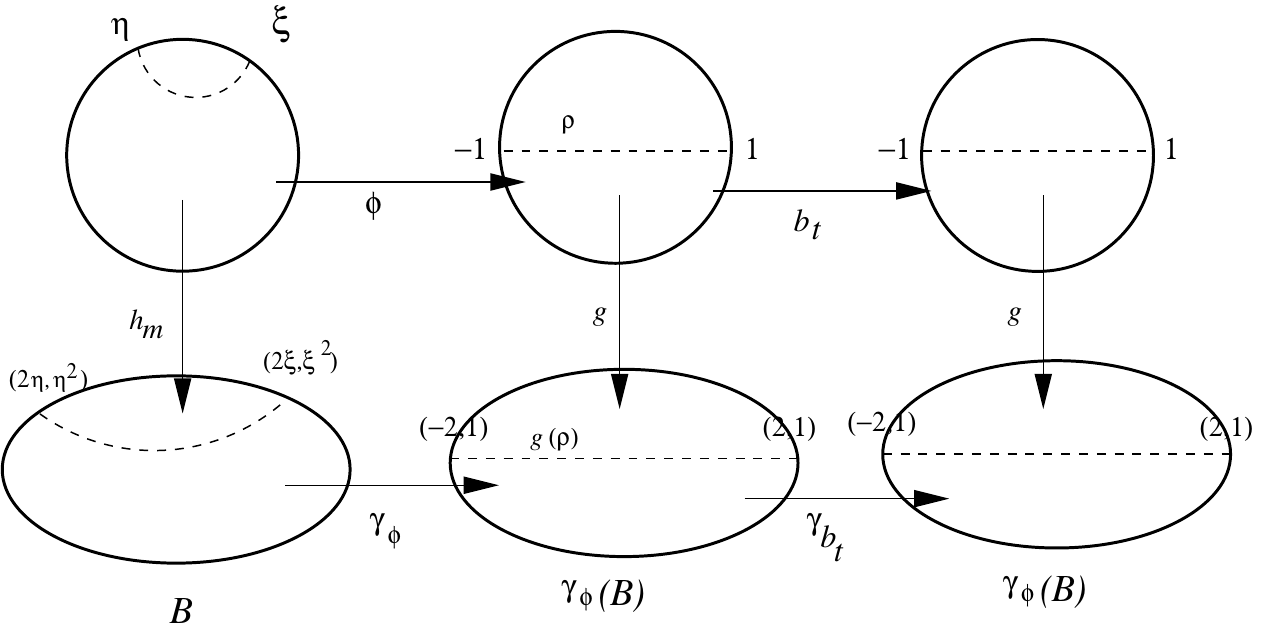}
\end{center}

Choose $\phi \in \aut \d$ satisfying $\phi(\xi) = 1$ and $\phi(\eta) =-1$. Since $\phi \circ m \circ \phi^{-1}(1) =1$ and $\phi \circ m \circ \phi^{-1} (-1) = -1$, it follows that there exists $r \in (-1,1)$ such that $\phi \circ m \circ \phi^{-1} = b_r$. Furthermore, $\gamma_\phi (B)$ is parametrized by the function $g$ defined by
\be\label{sharp.50}
g(z) =(z+b_r(z), zb_r(z)),\qquad z \in \d.
\ee

Since $b_r(-1) =-1$, $b_r(1) =1$, and $b_r|[-1,1]$ is strictly increasing, it follows from the Intermediate Value Theorem that there exists a unique $\rho \in (-1,1)$ such that $b_r(\rho) = -\rho$. By Lemma \ref{hyp.lem.8}, $b_r'(\rho) =1$, and we deduce from equation \eqref{sharp.50}  that
\[
g(\rho) = (0,-\rho^2)\ \ \text{ and }\ \ g'(\rho) = (2,0).
\]
Since $g$ parametrizes $\ga_\ph(B)$, $g'(\rho)\in T_{g(\rho)}\ga_\ph(B)$.
Therefore, by equation \eqref{sharpspecial}, $g(\rho)$ is a sharp point in $\gamma_\phi(B)$.

Now fix $t\in(-1,1)$.  Since $b_t \circ b_r = b_r \circ b_t$, 
\begin{align*}
\gamma_{b_t} \circ g(z)&=\gamma_{b_t}(z+b_r(z),zb_r(z))\\ 
&=(b_t(z) +b_t\circ b_r(z),b_t(z) b_t\circ b_r(z))\\ 
&=(b_t(z) +b_r\circ b_t(z),b_t(z) b_r\circ b_t(z)) =g\circ b_t(z),
\end{align*}
for all $z \in \d$. Hence
$\gamma_{b_t}\circ g(\rho) = g\circ b_t(\rho)$, and the tangent space
\[
T_{g\circ b_t(\rho) }\gamma_\phi(B) = \gamma_{b_t}'\circ g(\rho)\
T_{g(\rho)}\gamma_\phi(B).
\]
Therefore, as $g(\rho)$ is a sharp point in $\gamma\circ\phi(B)$,
\begin{align*}
T_{g\circ b_t(\rho)}\gamma\circ\phi(B) &= \gamma_{b_t}' \circ g(\rho)\ T_{g(\rho)}\gamma_\phi(B)\\ 
&= \gamma_{b_t}'\circ g(\rho) \ \sharp_{g(\rho)}=\sharp_{g\circ b_t(\rho)},
\end{align*}
that is,
\be\label{sharp.60}
g\circ b_t( \rho) \in \sharp(\gamma_\phi(B)).
\ee
Since equation \eqref{sharp.60} holds for all $t \in (-1,1)$ and 
\[
 \set{b_t(\rho)}{t\in(-1,1)} = (-1,1) = C(\{-1,1\}),
\] 
it follows that $g((-1,1)) \subseteq \sharp(\gamma_\phi(B))$, which is to say that
\[
\Xi_{\gamma_\phi(B)}  \subseteq \sharp(\gamma_\phi(B)).
\]
As $\gamma_\phi$ is an automorphism of $G$, $\Xi_B \subseteq \sharp(B)$, as was to be proved.
\end{proof}

\begin{defin}\label{dis.def.20}
 A {\em simple real geodesic} in $G$
\index{geodesic!simple real}
 is a curve of the form $h(C)$ where $h:\d \to G$ is a hyperbolic isometry and $C$ is a real geodesic in $\d$.

A \emph{distinguished geodesic in $G$}
\index{geodesic!distinguished}
 is a simple real geodesic in $G$ whose endpoints lie in the edge ${\mathrm E}$ of the distinguished boundary of $G$.
\end{defin}

Thus a simple real geodesic in $G$ is a curve $C$ in $G$ such that, for any pair of points $\la,\mu \in C$, the segment of $C$ joining $\la$ and $\mu$ achieves the minimum hyperbolic length of any curve in $G$ joining $\la$ and $\mu$. Note that not all real geodesics are simple.

Two types of distinguished geodesic
 $\Xi$ in $G$ are as follows. We say $\Xi$ is a \emph{royal distinguished} geodesic
\index{geodesic!royal distinguished}
 if $\Xi =h(C)$ where $h:\d \to \royal$ is a hyperbolic isometry and $C$ is a real hyperbolic geodesic in $\d$. We say $\Xi$ is a \emph{purely balanced distinguished geodesic}
\index{geodesic!purely balanced distinguished}
 if $\Xi = \Xi_B$ for some purely balanced geodesic $B$.
\begin{prop}\label{dis.prop.50}
$\Xi$ is a distinguished geodesic in $G$ if and only if it is either a royal or a purely balanced distinguished geodesic.
\end{prop}
\begin{lem}\label{dis.lem.20}
If $\lambda \in G\setminus \royal$, then there exists a purely balanced geodesic $B$ such that $\lambda \in \Xi_B$.
\end{lem}
\begin{proof}
Let $B$ be the geodesic in $G$ such that $\lambda \in B$ and $TB_\lambda = \sharp_\lambda$.  By Corollary \ref{sharp.cor.20}, $B$ is purely balanced. Also, since $\lambda \in \sharp(B)$, Proposition \ref{dis.prop.40} implies that $\lambda \in \Xi_B$.
\end{proof}
\begin{lem}\label{dis.lem.30}
If $B_1$ and $B_2$ are distinct purely balanced geodesics in $G$, then $\Xi_{B_1} \cap \Xi_{B_2} = \varnothing$.
\end{lem}
\begin{proof}
If $\lambda \in \Xi_{B_1} \cap \Xi_{B_2}$, then
\[
T_\la B_1= \sharp_\lambda = T_\la B_2.
\]
Hence, Theorem \ref{pre.thm.10} implies that $B_1=B_2$.
\end{proof}
The above two lemmas have the following immediate consequence.
\begin{thm}\label{dis.thm.20}
The purely balanced distinguished geodesics foliate $G\setminus \royal$.
\end{thm}
As a consequence of Theorem \ref{dis.thm.20} we may seek coordinates for $G \setminus \royal$ based on a parametrization of the purely balanced distinguished geodesics in $G$. Let
\[
X = \set{\{\xi_1,\xi_2\}}{\xi_1,\xi_2 \in E, \xi_1 \not=\xi_2}
\]
Since ${\mathrm E}$ is a circle, topologically, $X$  is homeomorphic to ${\mathrm M} \setminus {\mathrm E}$, a M\"obius band without boundary. For each point $\xi=\{\xi_1,\xi_2\} \in X$ we may construct a smooth parametrization
 \[
 t \mapsto B_{\xi,t},\qquad t\in (0,1),
 \]
of the family of purely balanced geodesics $B$ such that $B \cap {\mathrm E} = \xi$. Finally, for each $\xi \in X$ and each $t \in (0,1)$ we may smoothly parametrize $\Xi_{\xi,t}$ using $(0,1)$. These parametrizations lead to the following result.
\begin{thm}\label{dis.thm.30}
$G \setminus \royal$ is naturally homeomorphic to $X \times (0,1) \times (0,1) $ via the parametrizations described above. 
\end{thm}

%% file: closestpoint.tex

\section{The closest point property in $G$}

In \cite[Lemma 9.8]{aly2016} it was shown that if $F$ is a flat geodesic in $G$ and $\mu \in G \setminus F$, then
\be\label{cpp.10}
\inf_{\lambda\in F} d_G(\lambda,\mu) \mbox{ is attained at a point }  \lambda_0 \in F.
\ee
Further\footnote{Corollary \ref{cpp.cor.10} below asserts that $\lambda_0$ is unique.}, it was shown that if $\lambda_0$ and $\mu$ satisfy equation \eqref{cpp.10}, and $D$ denotes the geodesic passing through $\mu$ and $\lambda_0$, then
\be\label{cpp.20}
\mu \in \royal \implies D=\royal\ \ \text{ and }\ \ \mu\not\in \royal
\implies D \text{ is purely balanced.}
\ee
This result prompts the following definition.
\begin{defin}\label{ccp.def.10}
Let $F$ be a flat geodesic in $G$ and let $\mu\in G$. We say that $\lambda_0$ is \emph{a closest point in $F$ to $\mu$} if 
\[
d_G(\lambda_0,\mu)=\inf_{\lambda\in F} d_G(\lambda,\mu).
\]
If $D$ is a geodesic in $G$, we say that \emph{$D$ has the closest point property with respect to $F$} if $F$ meets $D$ in a point $\lambda_0$ and
for some $\mu \in D \setminus \{\lambda_0\}$,
\be\label{cpp.30}
\lambda_0 \text{ is a closest point in } F \text{ to } \mu .
\ee
\end{defin}
 Theorem \ref{cpp.thm.20} below implies that if $D$ has the closest point property with respect to $F$, then, in fact, statement \eqref{cpp.30} holds for every $\mu\in D$.
\subsection{Critical pairs}
\begin{defin}\label{cpp.def.20}
Let  $\lambda_0,\mu \in G$. We say that the pair $(\lambda_0,\mu)$ is a \emph{critical pair} if for every  admissible direction $v\in \flat_{\lambda_0}$,
\be\label{cpp.40}
\frac{d}{dt}\  d_G(\lambda_0 +tv,\mu)\  \big|_{t=0} \ge 0.
\ee
\end{defin}
\begin{lem}\label{cpp.lem.10}
Let $F$ be a flat geodesic in $G$ and let $\mu \in G \setminus F$. If $\lambda_0$ is a closest point in $F$ to $\mu$, then $(\lambda_0,\mu)$ is a critical pair.
\end{lem}
If $\lambda_0,\mu \in G$ and the geodesic $D$ passing through $\lambda_0$ and $\mu$ is purely balanced, then we may let $D=B_m$, where $m$ is hyperbolic with fixed points at $\eta_1$ and $\eta_2$. It then follows, by Lemma \ref{hyp.lem.30}, that
\[
d_G(\lambda_0,\mu) =\max_{\omega \in \{\bar\eta_1,\bar\eta_2\}}\
d(\Phi_\omega(\lambda_0),\Phi_\omega(\mu)).
\]
Consequently, if $v\in \flat_{\lambda_0}$, the following formula holds:
\be\label{cpp.50}
\frac{d}{dt}\ d_G(\lambda_0+tv,\mu) \big|_{t=0}=
\max_{\omega \in \{\bar\eta_1,\bar\eta_2\}}\
\frac{d}{dt}\ d(\Phi_\omega(\lambda_0+tv),\Phi_\omega(\mu))\ \big|_{t=0}.
\ee
\begin{lem}\label{cpp.lem.30}
Assume that $\lambda_0 =(0,p_0)$ is a point in $F_0$, the flat geodesic passing through the origin, and let $\mu \in G\setminus F_0$. Assume that the geodesic $D$ passing through $\lambda_0$ and $\mu$ is purely balanced and let $D=B_m$ where $m$ is hyperbolic with fixed points $\eta_1$ and $\eta_2$. $(\lambda_0,\mu)$ is a critical pair if and only if, for all $\xi \in \c$,
\be\label{cpp.60}
\mbox{ either } \quad \frac{d}{dt}\  d(p_0+t \xi, w_1)\Big|_{t=0}\ge 0 \quad \mbox{ or } \quad
\frac{d}{dt}\  d(p_0+t \xi, w_2)\Big|_{t=0}\ge 0,
\ee
where $w_1,w_2 \in \d$ are defined by
\be\label{cpp.70}
w_1=\eta_1\Phi_{\bar\eta_1}(\mu)\ \ \text{ and }\ \
w_2=\eta_2\Phi_{\bar\eta_2}(\mu).
\ee
\end{lem}
\begin{proof}
First notice that $v \in \flat_{(\lambda_0,\mu)}$ if and only if $v=(0,d)$ for some $d\in \c$. Consequently, using
Definition \ref{cpp.def.20} and \eqref{cpp.50}, we see that $(\lambda_0,\mu)$ is a critical pair if and only if for all $\xi\in \c$ there exists $\omega \in\{\bar\eta_1,\bar\eta_2\}$ such that 
\be\label{cpp.80}
\frac{d}{dt}\  d(\Phi_\omega(0,p_0+t\xi), \Phi_\omega(\mu))\Big|_{t=0}
\ge 0.
\ee
Note that, if $\omega \in \t$,
\begin{align*}
d(\Phi_\omega(0,p_0+t\xi), \Phi_\omega(\mu))&=
d(\omega(p_0+t\xi), \Phi_\omega(\mu))\\ 
&=d(p_0+t\xi, \bar\omega\Phi_\omega(\mu)).
\end{align*}
Hence statement \eqref{cpp.80} becomes
\[
\fa \xi\in \c \; \mbox{ there exists } \omega \in\{\bar\eta_1,\bar\eta_2\} \mbox{ such that }
\frac{d}{dt}\  d(p_0+t\xi, \bar\omega\Phi_\omega(\mu))\Big|_{t=0}
\ge 0,
\]
which is equivalent to statement \eqref{cpp.60}.
\end{proof}
A moment's thought reveals that statement \eqref{cpp.80} is equivalent to the assertion that $w_2$ is the hyperbolic reflection of $w_1$ about the point $p_0$. Thus, we obtain the following result.
\begin{lem}\label{cpp.lem.40}
Assume that $\lambda_0 =(0,p_0)$ is a point in  the flat geodesic $F_0$ passing through the origin, and let $\mu \in G\setminus F_0$. Assume that the geodesic $D$ passing through $\lambda_0$ and $\mu$ is purely balanced and let $D=B_m$ where $m$ is hyperbolic with fixed points $\eta_1$ and $\eta_2$.
Let $w_1=\eta_1\Phi_{\bar\eta_1}(\mu)$ and $w_2=\eta_2\Phi_{\bar\eta_2}(\mu)$.
 Then $(\lambda_0,\mu)$ is a critical pair if and only if $w_2$ is the hyperbolic reflection of $w_1$ about $p_0$. 
\end{lem}
\subsection{Critical pairs and orthogonality}
\begin{lem}\label{cpp.lem.50}
Let $p_0\in \d\setminus\{0\}$. A geodesic $D$ is purely balanced and meets $F_0$ at the point $(0,p_0)$ if and only if there exists $\alpha \in \d \setminus \{0\}$ and hyperbolic $m\in \aut \d$ such that  $D=B_m$, $m(\alpha)=-\alpha$, and $\alpha^2=-p_0$.
\end{lem}
\begin{proof}
 By Proposition \ref{hyp.prop.5}, $m$ is hyperbolic if and only if  $D=B_m$ is purely balanced. In addition, $m(\alpha) =-\alpha$ and $\alpha^2=-p_0$ if and only if $h_m(\alpha) =(0,p_0)$.
\end{proof}
\begin{lem}\label{cpp.lem.60}
Let $m\in \aut \d$ be hyperbolic and assume that $\alpha \in  \d \setminus \{0\}$ where $m(\alpha)=-\alpha$, and $\alpha^2=-p_0$.
There exists $\tau \in \t$ such that\footnote{By the Chain Rule, $\tau=m'(\alpha)$.} $m=b_\alpha \circ r_\tau \circ b_\alpha$.
 Furthermore, if we let $\eta_1,\eta_2$ denote the fixed points of $m$ and define
\[
\sigma=\eta_1+\eta_2\ \ \text{ and }\ \ \pi=\eta_1\eta_2,
\]
then
\be\label{cpp.90}
\sigma =\frac{1-|\alpha|^2}{\bar\alpha}\ \frac{1-\tau}{1+\tau}
\ \ \text{ and }\ \ \pi= -\frac{\alpha}{\bar\alpha}.
\ee
\end{lem}
\begin{proof}
That the general $m \in \aut \d$ satisfying $m(\alpha)=-\alpha$ has the form $m=b_\alpha \circ r_\tau \circ b_\alpha$ follows by Schur reduction.

To prove equations \eqref{cpp.90} it suffices to show that the equation $m(x)=x$ is equivalent to the equation $x^2-\sigma x +\pi=0$.
\begin{align*}
m(x)=x \iff\ \ & (b_\alpha \circ r_\tau \circ b_\alpha) (x) = x\\ 
\iff\ \ &\tau b_\alpha(x)=b_{-\alpha}(x)\\ 
\iff\ \ & \tau(x-\alpha)(1+\bar\alpha x)=(x+\alpha)(1-\bar\alpha x)\\ 
\iff\ \ &\bar\alpha(1+\tau) x^2 -(1-|\alpha|^2)(1-\tau)x -\alpha(1+\tau)=0\\ 
\iff \ \ &x^2-\frac{1-|\alpha|^2}{\bar\alpha}\ \frac{1-\tau}{1+\tau}x -\frac{\alpha}{\bar\alpha}=0\\
\iff\ \ & x^2-\sigma x +\pi =0.
\end{align*}
\end{proof}
We remark that equations \eqref{cpp.90}  imply that $\tau=1$ if and only if 
$\sigma=0$, a fact which also follows from Lemma \ref{hyp.lem.7a} and Condition (ii) in Lemma \ref{hyp.lem.8}. Also observe that the formulas in  \eqref{cpp.90} imply the formulas
\be\label{cpp.95}
\frac{1+\tau}{1-\tau}\sigma=\frac{1-|\alpha|^2}{\bar\alpha}
\ \ \text{ and, since }\ \tau \in \t, \ \
\frac{1+\tau}{1-\tau}\bar\sigma=-\frac{1-|\alpha|^2}{\alpha}.
\ee
\begin{lem}\label{cpp.lem.70}
Assume that $\lambda_0 =(0,p_0)$ is a point in $F_0$, the flat geodesic passing through the origin, and let $\mu \in G\setminus F_0$. Assume that the geodesic $D$ passing through $\lambda_0$ and $\mu$ is purely balanced and let $D=B_m$ where $m$ is hyperbolic with fixed points $\eta_1$ and $\eta_2$. If $(\lambda_0,\mu)$ is a critical pair, then $D$ is orthogonal to $F_0$.
\end{lem}
\begin{proof}
By Lemma \ref{cpp.lem.50}, there exists $\alpha \in \d$ such that
\be\label{cpp.100}
m(\alpha) =-\alpha,\ \ \text{ and }\ \ \alpha^2=p_0,
\ee
equations that imply
\be\label{cpp.110}
h_m(\alpha) =(0,p_0).
\ee
As also $\mu \in B_m$, there exists $\beta \in \d$ such that
\be\label{cpp.120}
h_m(\beta) =\mu.
\ee

Lemma \ref{hyp.lem.30} guarantees that if we define $m_1$ and $m_2$ by
\be\label{cpp.130}
m_1 =(\eta_1\Phi_{\bar\eta_1}) \circ h_m\ \ \text{ and }\ \
m_1 =(\eta_2\Phi_{\bar\eta_2}) \circ h_m
\ee
then $m_1,m_2 \in \aut \d$. Observe that
\[
m_1(\alpha)=\eta_1\Phi_{\bar\eta_1}(0,p_0)=p_0,
\]
and using equations $\eqref{cpp.70}$, that
\[
m_1(\beta) = \eta_1\Phi_{\bar\eta_1}(\mu) =w_1.
\]
Similarly,
\[
m_2(\alpha) =p_0\ \ \text{ and }\ \ m_2(\beta)=w_2.
\]
Consequently, if we set $\phi=m_2\circ m_1^{-1}$,
\[
\phi(p_0) = p_0\ \ \text{ and }\ \ \phi(w_1)=w_2.
\]
As Lemma \ref{cpp.lem.40} guarantees that $w_2$ is the hyperbolic reflection of $w_1$ about $p_0$, it follows that $\phi$ is hyperbolic reflection about $p_0$. Therefore, as $\phi'(p_0)=-1$,
\be\label{cpp.140}
m_1'(\alpha)+m_2'(\alpha)=0
\ee
Noting that, for $\eta\in \t$,
\[
\nabla (\eta\Phi_{\bar\eta})(s,p)=\left(\frac{-2\eta+2\bar\eta p}{(2-\bar\eta s)^2},\frac{2}{2-\bar\eta s}\right),
\]
we see that
\begin{align*}
m_1'(\alpha)&=\nabla (\eta_1\Phi_{\bar\eta_1})(0,p_0)\ \cdot \ h_m'(\alpha)\\ 
&=\big(-\tfrac12\eta_1+\tfrac12 \bar\eta_1 p_0,1\big) \cdot \big(1+m'(\alpha),\alpha(m'(\alpha)-1)\big)\\ 
&=(1+m'(\alpha))(-\tfrac12\eta_1+\tfrac12 \bar\eta_1 p_0)+
\alpha(m'(\alpha)-1).
\end{align*}
Likewise
\[
m_2'(\alpha) = (1+m'(\alpha))(-\tfrac12\eta_2+\tfrac12 \bar\eta_2 p_0)+
\alpha(m'(\alpha)-1).
\]
Therefore, by equation \eqref{cpp.140} and Lemma \ref{cpp.lem.60},
\begin{align*}
0&\ =\ \ m_1'(\alpha) +m_2'(\alpha)\\ 
&\ =\ \ (1+m'(\alpha))(-\tfrac12\eta_1+\tfrac12 \bar\eta_1 p_0)+
\alpha(m'(\alpha)-1)\\
&\ \ \ \ +(1+m'(\alpha))(-\tfrac12\eta_2+\tfrac12 \bar\eta_2 p_0)+
\alpha(m'(\alpha)-1)\\ 
&\ =\ \ (1+\tau)(-\tfrac12 \sigma -\tfrac12 \bar\sigma \alpha^2)+2\alpha(\tau-1)\qquad \mbox{ by Lemma  \ref{cpp.lem.60} }\\ 
&\ =\ \ -\frac{1-\tau}{2}\Big(\frac{1+\tau}{1-\tau}\sigma+
\frac{1+\tau}{1-\tau}\bar\sigma\alpha^2+4\alpha\Big)\\ 
&\ =\ \  -\frac{1-\tau}{2}\Big(\frac{1-|\alpha|^2}{\bar\alpha}
-\frac{1-|\alpha|^2}{\alpha}\alpha^2+4\alpha\Big)\quad \mbox{ by equations \eqref{cpp.95} }\\ 
&\ =\ \ -\frac{1-\tau}{2\bar\alpha}\big(1-|\alpha|^2\big)^2.
\end{align*}
Consequently, as $m'(\alpha) =\tau$, $m'(\alpha)=1$ and  Lemma \ref{hyp.lem.8} implies that $m$ is irrotational. By Proposition \ref{sharp.prop.10}, $D$ is orthogonal to $F_0$.
\end{proof}
\subsection{The closest point property and orthogonality}
\begin{thm}\label{cpp.thm.20}
Let $F$ be a flat geodesic in $G$. Assume that $D$ is a geodesic in $G$ that meets $F$ at the point $\lambda_0$. The following are equivalent.
\begin{enumerate}[\rm (i)]
\item For every $\mu \in D$ $\lambda_0$ is a closest point in $F$ to $\mu$;
    \item there exists $\mu \in D$ such that $\lambda_0$ is a closest point in $F$ to $\mu$;
            \item $D$ is orthogonal to $F$.
\end{enumerate}
\end{thm}
\begin{proof}
Clearly, (i) implies (ii). If (ii) holds, then \eqref{cpp.20} implies that $D$ is purely balanced and
Lemma \ref{cpp.lem.10} implies that $(\lambda_0,\mu)$ is a critical pair. Therefore, Lemma \ref{cpp.lem.70} implies that (iii) holds. 

Now suppose that (iii) holds. Fix $\mu\in \d$. Choose $\lambda_1\in F$ such that $\lambda_1$ is a closest point in $F$ to $\mu$.
By (ii) implies (iii), if $D_1$ is the geodesic passing through $\lambda_1$ and $\mu$, then $D_1$ is orthogonal to $F$. As $D$ and $D_1$ are both orthogonal to $F$ and $\mu\in D\cap D_1$, Lemma \ref{sharp.lem.90} implies that  $\lambda_1=\lambda_0$. Therefore, $\lambda_0$ is a closest point in $F$ to $\mu$.
\end{proof}
\begin{corollary}\label{cpp.cor.10}
Assume that $F$ is a flat geodesic in $G$ and $\mu\in G$. If $\lambda_1$ and $\lambda_2$ are both a closest point in $F$ to $\mu$, then $\lambda_1=\lambda_2$.
\end{corollary}
\begin{proof}
Suppose
\[
\inf_{\lambda\in F} d_G(\lambda,\mu)
\]
is attained at two points $\lambda_1$ and $\lambda_2$ in $F$. If we let $D_1$ be the geodesic passing through $\lambda_1$ and $\mu$ and let $D_2$ be the geodesic passing through $\lambda_2$ and $\mu$, then the equivalence of Conditions (ii) and (iii) in Theorem \ref{cpp.thm.20} implies that $D_1$ and $D_2$ are both orthogonal to $F$. Therefore, as $\mu \in D_1\cap D_2$, Lemma \ref{sharp.lem.90} implies that $D_1=D_2$. Hence, $\lambda_1=\lambda_2$.
\end{proof}
The following corollary of Theorem \ref{cpp.thm.20} is equivalent to Theorem \ref{int.thm.20} from the introduction of the paper.
\begin{corollary}\label{cpp.cor.20}
Let $F$ be a flat geodesic in $G$ and let $D$ be a geodesic in $G$. Then $D$ has the closest point property with respect to $F$ if and only if $D$ is orthogonal to $F$.
\end{corollary}